\documentclass[11pt,reqno]{amsart}
\usepackage{amsmath,amsthm,amssymb,latexsym}

\newcommand{\psp}{\vspace{0.28cm}} 
\newcommand{\pse}{\vspace{0.4cm}} 

\newtheorem{thm}[equation]{Theorem}
\newtheorem{cor}[equation]{Corollary}
\newtheorem{lem}[equation]{Lemma}
\newtheorem{prop}[equation]{Proposition}
\newtheorem{exam}[equation]{Example}

\newtheorem{rem}[equation]{Remark}
\newtheorem{rems}[equation]{Remarks}

\theoremstyle{definition}
\newtheorem{defn}[equation]{Definition}

\numberwithin{equation}{section}

\newcommand{\ot}{\otimes}
\newcommand{\wot}{\widetilde \otimes}

\newcommand{\la}{\langle}
\newcommand{\ra}{\rangle}
\newcommand{\pa}{\partial}

\newcommand{\al}{\alpha}
\newcommand{\be}{\beta}
\newcommand{\Ga}{\Gamma}
\newcommand{\ga}{\gamma}

\newcommand{\De}{\Delta}

\newcommand{\lam}{\lambda}
\newcommand{\La}{\Lambda}
\newcommand{\si}{\sigma}
\newcommand{\sig}{\si}
\newcommand{\Th}{\theta}
\newcommand{\om}{\omega}

\newcommand{\Z}{\mathbb Z}
\newcommand{\ZZ}{\Z}
\newcommand{\CC}{\mathbb C}
\newcommand{\F}{\mathbb F}
\newcommand{\K}{\mathbb K}
\newcommand{\NN}{\mathbb N}

\newcommand{\g}{{\mathfrak g}}
\newcommand{\GG}{{\mathfrak T}}
\newcommand{\FF}{{\mathfrak F}}
\newcommand{\h}{{\mathfrak h}}

\newcommand{\scl}{{\mathfrak L}}
\newcommand{\KK}{{\mathcal K}}
\newcommand{\LL}{{\mathcal L}}
\newcommand{\MM}{{\mathfrak M}}
\newcommand{\RR}{{\mathfrak R}}

\newcommand{\fru}{{\mathfrak u}}

\newcommand{\sca}{{\mathcal A}}
\newcommand{\scA}{\sca}
 \newcommand{\scB}{{\mathcal B}}
\newcommand{\scC}{{\mathcal C}}
\newcommand{\scd}{{\mathcal D}}\newcommand{\scD}{{\mathcal D}}
\newcommand{\scE}{{\mathcal E}}
\newcommand{\scH}{{\mathcal H}}
\newcommand{\scI}{{\mathcal I}}
\newcommand{\scJ}{{\mathcal J}}
\newcommand{\scK}{{\mathcal K}}
\newcommand{\scL}{{\mathcal L}}

\newcommand{\scM}{{\mathcal M}}
\newcommand{\scN}{{\mathcal N}}
\newcommand{\scQ}{{\mathcal Q}}

\newcommand{\scV}{{\mathcal V}}

\newcommand{\Ann}{\mathord{\rm Ann}}
\newcommand{\ad}{\mathop{\rm ad}}
\newcommand{\Aut}{\mathord{\rm Aut}}
\newcommand{\rmC}{{\rm Cent}}
\newcommand{\cent}{\rmC}
\newcommand{\grcent}{\mathord{\rm grCent}}
\newcommand{\Der}{\mathord{\rm Der}}
\newcommand{\End}{\mathord{\rm End}}
\newcommand{\grEnd}{\mathord{\rm grEnd}}
\newcommand{\ev}{\mathord{\tt ev}}
\newcommand{\Hom}{\mathord{\rm Hom}}
\newcommand{\Id}{\mathop{\rm id}}
\newcommand{\im}{\mathop{\rm im}\,}
\newcommand{\Ker}{\ker}
\newcommand{\mult}{\mathop{\rm Mult}}
\newcommand{\rad}{\mathord{\rm rad}\,}
\newcommand{\SDer}{\mathord{\rm grSDer}}
\newcommand{\supp}{\mathord{\rm supp}\,}

\newcommand{\Proof}{\proof}
\newcommand{\ms}{\medskip}

\begin{document}

\title[]{\large THE CENTROID OF  EXTENDED AFFINE AND \\
ROOT GRADED LIE ALGEBRAS}

\author[]{Georgia Benkart$^{a,1}$}
\thanks{$^{1}$Supported in part by National Science Foundation Grant
\#{}DMS--0245082.} 
\author[]{Erhard Neher$^{b,2,*}$}
\thanks{$^{2}$Supported in part by Natural Sciences and Engineering
Research Council of Canada Discovery Grant \#{}8836--2001. 
\newline{$^{*}$Corresponding author}
\newline{{\it E-mail addresses:} benkart@math.wisc.edu (G. Benkart), 
    neher@uottawa.ca (E. Neher)} 
\newline 2000 Mathematical Subject
Classification: Primary 17B40; Secondary 17A36} 
\date{\today}
\maketitle

\address{\small $^A$\it Department of Mathematics, University of 
Wisconsin, Madison, WI 53706-1388 USA} 

\address{\small $^B$\it Department of Mathematics and Statistics, 
University of Ottawa, Ottawa, Ontario, K1N 6N5, Canada } 

%
%
%
%

\begin{abstract}
We develop general results on centroids of Lie algebras and apply 
them to determine the centroid of extended affine Lie algebras, 
loop-like  and Kac-Moody Lie algebras, and Lie algebras graded by 
finite root systems.
\end{abstract}

\section{Introduction}
 
 Our main focus will be    
on centroids of Lie algebras.
 When $\scL$ is a Lie algebra,  the centroid  $\cent(\scL)$ is just the space of $\scL$-module homomorphisms $\chi$ on $\scL$:  \  $\chi([x,y]) = [x,\chi(y)]$ for all $x,y \in \scL$,  (viewing $\scL$ as an $\scL$-module under the adjoint action).   Our interest in the centroid stems from investigations of extended affine Lie algebras   (see \cite{AABGP}).
These Lie algebras are natural generalizations of the affine and 
toroidal Lie algebras,  which have played such a critical role in 
many different areas of mathematics and physics.   Their root 
systems (the so-called extended affine root systems) feature 
prominently in the work of Saito (\cite{S1}, \cite{S2}) and Slodowy 
\cite{Sl} on singularities. In the classification of the extended 
affine Lie algebras, elements of the centroid are essential in 
constructing the portion of the algebra that lies outside of the 
core (see for example, \cite{N2}). This is the part of the extended 
affine Lie algebra $\scE$ that is nondegenerately paired with the 
centre under the invariant bilinear form on $\scE$.  Thus, results 
on the centroid are a key ingredient in the classification of 
extended affine Lie algebras. 

The centroid also plays an important role in understanding forms of 
an algebra: All scalar extensions of a simple algebra remain simple 
if and only if its centroid just consists of the scalars in the base 
field. In particular, for  finite-dimensional  simple associative 
algebras,  the centroid is critical in investigating Brauer groups 
and division algebras. Another area where the centroid occurs 
naturally is in the study of derivations of an algebra. If  $\chi 
\in \cent(\scA)$ and $\pa$ is a derivation of $\scA$, then $\chi 
\pa$ is also a derivation of $\scA$, so centroidal transformations 
can be used to construct derivations of an algebra. \psp 

We will develop general results, extending some earlier work of 
 other authors, and then apply them to determine the centroid of several families of 
(mostly infinite-dimensional) Lie algebras: \  extended affine Lie 
algebras, loop algebras, Lie algebras graded by finite root systems, 
and Kac-Moody Lie algebras,  which are distinguished because of 
their substantial applications in a diverse array of subjects.  \pse

\section{Centroids of algebras}\label{gen} 

\subsection{Some general results}\label{secgen}

We begin with a little background on centroids for arbitrary (not 
necessarily Lie, associative, etc.) algebras $\scA$. Proofs of the 
results quoted here can be found for example in \cite[X.1]{jake} for 
finite-dimensional algebras and in \cite[II, 1.6, 1.7]{Mc} in 
general. 

\psp It is natural and important for our approach to centroids to 
consider algebras over a unital commutative associative ring,  (for 
example,  a perfect algebra $\scA$ over its centroid,  which may not 
be a field unless $\scA$ is simple).  Thus, let $\scA$ be an 
arbitrary algebra over a unital commutative associative ring $\K$. 
The {\em centroid} of $\scA$ is the space of $\K$-linear 
transformations on $\scA$ given by $\cent(\scA) = \{\chi \in 
\End_\K(\scA) \mid \chi(ab) = a \chi(b) = \chi(a)b$ for all $a,b \in 
\scA \}$. We will write $\cent_\K(\scA)$ for $\cent(\scA)$ if it is 
important to emphasize the dependence on $\K$. Clearly, 
$\cent(\scA)$ is simply the centralizer algebra of the {\it 
multiplication algebra}  $\mult(\scA)= \mult_\K(\scA)$, the unital 
subalgebra of $\End_\K(\scA)$ generated by the left and 
multiplication operators of $\scA$. The centroid is always a 
subalgebra of the associative algebra $\End_\K(\scA)$. \psp

 For $a,b,c\in \scA$,  the {\it associator} is defined as 
$(a,b,c)=(ab)c-a(bc)$.  The {\it centre} of $\scA$ consists of all 
$z\in \scA$ satisfying $za=az$ and $(a,b,z)=(a,z,b)=(z,a,b)=0$ for 
all $a,b\in \scA$. The centre is always a commutative associative 
subalgebra of $\scA$. Moreover, if $\scA$ has an identity element 
$1$,  then $\cent(\scA) \to \scA, \  \chi \mapsto \chi(1)$,  is an 
algebra isomorphism between the centroid and the centre of $\scA$. 
We denote by $\scA^{(1)}$ the $\K$-span of all products $ab$ for 
$a,b \in \scA$. If $\scA$ is {\it perfect}, (i.e., $\scA$ equals 
$\scA^{(1)}:=\scA\scA$), then the centroid is necessarily 
commutative, as $\chi \psi(ab) = \chi(\psi(a)b) = \psi(a)\chi(b) = 
\psi(a \chi(b)) = \psi \chi(ab)$ holds for all $a,b \in \scA$, 
$\chi,\psi \in \cent(\scA)$. If $\cent(A)$ is commutative, we may 
regard $\scA$ as an algebra over its centroid,  $\chi a = a \chi$ 
and $\chi(ab) = (\chi a)b = a (\chi b)$ for all $a,b \in \scA$, 
$\chi \in \cent(\scA)$.   If $\scA$ is {\it prime} in the sense that 
$\scA$ has no nonzero ideals $I$, $J$ with $I J=0$, then 
$\cent(\scA)$ is an integral domain and $\scA$ is a torsion-free 
$\cent(A)$-module.  If $\scA$ is {\it simple}, i.e., $\scA^{(1)}\ne 
0$,  and the only ideals of $\scA$ are  $\scA$ and $0$, then 
$\cent(\scA)$ must be a field by Schur's Lemma.  When the centroid 
of an algebra coincides with the base ring $\K$ (more precisely, it 
equals $\K\Id$), the algebra is said to be {\em central},  and in 
the special case of a simple algebra, it is said to be {\em central 
simple}.  Every simple algebra is central simple over its centroid.  
\psp

 For  any subset $\scB$ of $\scA$, the  {\it annihilator of 
$\scB$ in $\scA$}  is $\Ann_\scA(\scB)= \{ z\in \scA \mid 
z\scB=0=\scB z\}$.  Any $\K$-submodule of $\scA$ containing 
$\scA^{(1)}$ or contained in $\Ann_{\scA}(\scA)$ is an ideal of 
$\scA$.  If $\scA$ is a Lie algebra,  we follow the usual convention 
of denoting the product of $a,b\in \scA$ by $[a,b]$.  In this case, 
 $\Ann_{\scA}(\scB)$ is simply  the {\it centralizer} $\scC_\scA(\scB) = \{z \in \scA \mid [z,\scB] =0\}$ of $\scB$ in $\scA$. 
In particular, $\Ann_{\scA}(\scA)= Z(\scA)$, the usual centre of 
$\scA$, which coincides with the general definition of the centre as 
given above if $\frac{1}{2} \in \K$.  Let $\Der(\scA)$ denote the 
algebra $\Der_{\K}(\scA)$ of $\K$-linear derivations of an algebra 
$\scA$.   Then we have the following basic facts: \psp

\begin{lem}\label{easy}  Let $\scA$ be an algebra over 
a unital commutative associative ring $\K$ and let $\scB$ be a 
subset of $\scA$. 
\begin{itemize} 
\item[{\rm (a)}]  $\Ann_\scA (\scB)$ is invariant 
under  $\cent(\scA)$,  as is  any perfect  ideal of $\scA$. 
\smallbreak 
 
\item[{\rm (b)}]   For any $\cent(\scA)$-invariant ideal $\scB$ of $\scA$,  
the {\rm vanishing ideal} $${\scV}(B) :=  \{ \chi\in \cent(\scA) 
\mid \chi(\scB)=0 \}$$  
 is isomorphic to $\Hom_{\scA/\scB}(\scA/\scB, {\Ann}_\scA(\scB))$, which is the 
set of $\K$-linear maps $f : \scA/\scB \to \Ann_\scA(\scB)$ 
satisfying $f(xy) = f(x)y = xf(y)$ for all $x,y\in \scA/\scB$, where 
$\scA/\scB \, \times \, \Ann_\scA(\scB) \to \scA$ is defined by $(a+ 
\scB)z = az $ and similarly for $\Ann_\scA(\scB) \, \times \,  
\scA/\scB \to \scA$. In particular,  if $\cent(\scB) = \K \Id$,  
then  
\begin{equation}\label{smacent} \cent(\scA) = \K \Id  \oplus 
\scV(\scB).  
 \end{equation} 
\medbreak 

\item[{\rm (c)}] \label{cideal}  \vspace{-.5cm} \begin{eqnarray}
    \cent(\scA) \cap \Der (\scA)
  &=& \{ \psi \in {\End}_\K(\scA) \mid \scA^{(1)}  \subseteq \ker \psi,
     \,\im \psi \subseteq {\rm Ann}_{\scA}(\scA) \} \nonumber \\
      &=& \{ \psi\in \cent(\scA) \mid \scA^{(1)} \subseteq \ker \psi \}
            = \scV(\scA^{(1)})\nonumber  \\
      &=& \{\psi \in \Der(\scA) \mid \im \psi \subseteq{\Ann}_{\scA}(\scA) \}, \nonumber \\
     &\cong& \Hom_{\scA/{\scA}^{(1)}}(\scA/\scA^{(1)} \,,\, {\Ann}_{\scA}(\scA))\nonumber\\
     &\cong&
     \Hom_\K(\scA/\scA^{(1)} \,,\, {\Ann}_{\scA}(\scA)).\nonumber
  \end{eqnarray}

\item[{\rm (d)}] $\scA$ is indecomposable (cannot be written as the direct sum of
two nontrivial ideals, or equivalently, is an indecomposable 
$\mult(\scA)$-module)  if and only if $\cent(\scA)$ does not contain 
idempotents $\ne 0, \Id$.  \smallbreak

\item[{\rm (e)}]  Suppose  $\scA$ is an 
indecomposable $\mult(\scA)$-module of finite length $n$,  and 
denote by $\rad \cent(\scA)$ the Jacobson radical of $\cent(\scA)$. 
Then $\cent(\scA)$ is a local ring, i.e., $\cent(\scA)/\rad 
\cent(\scA)$ is a division ring (see for example,  {\rm 
\cite[Sec.~19]{L}}), 
 and $\big(\rad \cent(\scA)\big)^n=0$. Thus,  $\rad 
\cent(\scA)$ is nilpotent and coincides with the set of nilpotent 
transformations in $\cent(\scA)$.  In particular, if $\scA$ is a 
finite-dimensional indecomposable  algebra over a perfect field $\F$ 
then there exists a division algebra $\mathfrak D$ over $\F$ such 
that $\cent(\scA) = \allowbreak \mathfrak D \Id \,\oplus \,\rad 
\cent(\scA)$. \smallbreak

\item[{\rm (f)}]  If $\scA$ is perfect, every $\chi \in \cent(\scA)$ is symmetric with 
respect to any invariant form on $\scA$.\end{itemize}
\end{lem}

\Proof (a) -- (d) are straightforward. Part (d)  can be found in  
\cite[Sec.~1]{Me} for Lie algebras or in \cite[Lem.~1]{Po} for more 
general algebras, where a description of $\cent(\scA)$  for 
decomposable $\scA$ is also given. 

(e) Since $\cent(\scA)$ consists of  the $\mult(\scA)$-module 
endomorphisms of $\scA$, the first part of (e) follows from 
\cite[Thm.~19.17]{L}. Under the assumptions of the second part we 
know that $\cent(\scA)$ is a local $\F$-algebra. The claim then 
follows from Wedderburn's Principal Theorem. 

(f) Let $(\,|\,)$ be an {\em invariant} $\K$-bilinear form on 
$\scA$, so that $(ab\,|\, c) = (a\,|\,bc)$ for all $a,b,c \in \scA$,  
and let $\chi \in \cent(\scA)$.  Then  $(\chi(ab)\,|\,c) = 
(a\chi(b)\,|\,c) = (a\,|\,\chi(b)c) = (a\,|\,b\chi(c))\allowbreak  = 
(ab\,|\,\chi(c))$ for all $a,b,c \in \scA$. \qed \psp

\begin{rem}\label{remcoho}  {\rm  Let $\scL$ be a Lie 
algebra over a field $\F$. A derivation from $\LL$ to an 
$\LL$-module $\mathfrak M$ is an $\F$-linear map $\delta: \LL \to 
\mathfrak M$ such that $\delta \bigl([x,y]\bigr) = x.\delta(y) - 
y.\delta(x)$  for all $x,y \in \LL$. The space $\Der(\LL,\mathfrak 
M)$ of such derivations contains the inner derivations  $\hbox{\rm 
IDer}(\LL, \mathfrak M) = \{\delta_m \mid  m \in \mathfrak M\}$,  
where  $\delta_m(x) = x.m$ for all $x \in \LL$.  Then the {\it first 
cohomology group of $\LL$ with values in $\mathfrak M$} is the 
quotient ${\tt H}^1(\LL,\mathfrak M) := \Der(\LL,\mathfrak 
M)/\hbox{\rm IDer}(\LL, \mathfrak M)$ (see for example,  
\cite[Ch.~I, Sec.~3, Ex.~12]{bou:lieI} or \cite[Ch.~V.6]{jake}).  
Now for any Lie algebra $\LL$,  an ideal $\MM$ of $\LL$ is an 
$\LL$-module under the adjoint action. Examples of 
$\cent(\LL)$-invariant ideals are the centre $Z(\LL)$ and all the 
ideals in the derived series, lower (descending) central series, and 
ascending central series of $\LL$.  In particular  if $\MM=Z(\LL)$,  
we have  $\hbox{\rm IDer}(\LL, Z(\LL)) = 0$. Thus ${\tt 
H}^1(\LL,Z(\LL))\allowbreak= \Der(\LL, Z(\LL))$,  and by Lemma 
\ref{cideal},  we have a canonical identification
\begin{eqnarray}\hspace{-.3truein}{\tt H}^1(\LL,Z(\LL))&=&\left \{\psi \in {\End}_\F(\LL)\, \big | \, [\psi(\LL), \LL] = 0 = \psi(\LL^{(1)})\right \}   = {\scV}(\scL^{(1)}) \label{eq:vl}\\
       &\cong& \Hom_\F\big(\LL/\LL^{(1)},Z(\scL)\big) \nonumber  \end{eqnarray}
as $\cent(\LL)$-modules.}
\end{rem} \ms

\begin{exam}\label{small} {\rm 
For any Lie algebra $\scL$ over a field $\F$, 
$${\scV}(\scL^{(1)}) =\left \{\psi \in {\End}_\F(\LL)\, \big | \, [\psi(\LL), \LL] = 0 = \psi(\LL^{(1)})\right \},$$
as in \eqref{eq:vl}.  Thus, if  $Z(\scL) \ne 0$ and $\scL \ne 
\scL^{(1)}$,  we have  
$$ \F\Id \subsetneq \F\Id \oplus {\scV}(\scL^{(1)}) \subseteq \cent_\F(\scL).$$
So in order for $\scL$ to be central, a necessary condition is that 
$Z(\scL)=0$ or $\scL$ is perfect.  Later results  (Corollaries 
\ref{toralcor} and \ref{ealacent})  will treat various classes of 
Lie algebras for which $\F\Id\subsetneq \F\Id \oplus 
{\scV}(\scL^{(1)}) = \cent(\scL)$.        {\it Heisenberg Lie 
algebras}  provide easy examples of this phenomenon.  A Heisenberg 
Lie algebra $\scH$ has a basis $\{a_i,b_i \mid i \in \mathcal I\} 
\cup \{c\}$, such that $[a_i,b_j] = \delta_{i,j}c, \ [a_i,a_j] = 0 = 
[b_i,b_j]$, and $[\scH,c]=0$, where $\delta_{i,j}$ is the Kronecker 
delta.    By (2.2) applied to $\scB = \scH^{(1)} = Z(\scH)$, $$ 
\cent(\scH) = \F\Id \oplus {\scV}(\scH^{(1)}).$$ For nilpotent Lie 
algebras  (in particular,  for ${\scL} = {\scH}$),   
${\scV}(\scL^{(1)}) \neq 0$, but it may well be that $\F \Id \oplus 
{\scV}(\scL^{(1)})\allowbreak  \subsetneq \cent({\scL})$ for an 
arbitrary nilpotent Lie algebra.  (see  \cite[Prop.~2.7]{Me}).} 
\end{exam} \psp

\begin{rem} {\rm Indecomposable Lie algebras $\scL$  having a
small centroid, i.e. those for which $\cent(\scL) = \F \Id \oplus 
{\scV}(\scL^{(1)})$,  have been investigated by Melville \cite{Me} 
and Ponomar{\"e}v  \cite{Po} under certain assumptions  (e.g. in  
\cite{Me}, when  $\scL$ is finite-dimensional).}  
\end{rem} \psp

\begin{lem}\label{oblem}
Let $\pi: \scA \to \scB$ be an epimorphism of $\K$-algebras. For 
every $f\in \End_\K(\scA;\Ker \pi):= \{ g\in \End_\K(\scA)\mid 
g(\Ker \pi) \subseteq \Ker \pi\}$,  there exists a unique $\bar f 
\in \End_\K(\scB)$ satisfying $\pi \circ f = \bar f \circ \pi$.  
Moreover the following hold:

 \begin{itemize}
 \item[{\rm (a)}]  The map
 $$ \pi_{\End} : {\End}_\K(\scA;\Ker \pi) \to
\hbox {\rm End}_\K(\scB), \quad  f \mapsto \bar f
 $$
is a unital algebra homomorphism with the following properties:
\begin{eqnarray}
\label{pimult} \pi_{\End} (\mult(\scA)) &=& \mult(\scB),  \\
\pi_{\End}\big(\cent(\scA) \cap {\End}_\K(\scA;\Ker \pi)\big) 
&\subseteq& \cent(\scB).  \label{picent}
\end{eqnarray}
By restriction, there is an algebra homomorphism
\begin{equation}
\pi_{\cent} : \big(\cent(\scA) \cap {\hbox{\rm End}_\K(\scA;\Ker 
\pi)\big)} \to
 \cent (\scB), \quad  f \mapsto \bar f
\end{equation}
If $\ker \pi = \Ann_{\scA}(\scA)$,  every $\chi \in \cent(\scA)$ 
leaves $\Ker \pi$ invariant,  and hence $\pi_{\cent}$ is defined on 
all of $\cent(\scA)$.
\smallskip

 \item[{\rm (b)}]  Suppose $\scA$ is perfect and $\Ker\pi \subseteq
 \Ann_{\scA}(\scA)$.
Then
\begin{equation}
\pi_{\cent} : \big(\cent(\scA) \cap  \hbox{\rm End}_\K(\scA;\Ker 
\pi)) \to
 \cent (\scB), \quad  f \mapsto \bar f
\end{equation}
is injective. 
\smallskip
\item[{\rm (c)}]  If $\scA$ is perfect, $\Ann_{\scB}(\scB) = 0$ and $\Ker\pi
\subseteq \Ann_{\scA}(\scA)$, then $\pi_{\cent} : \cent(\scA) \to 
\cent(\scB)$  is an algebra monomorphism.
  \end{itemize}
\end{lem}

The main application of this lemma will be to Lie algebras. In this 
case, an epimorphism $\pi : \scA \to \scB$ with $\Ker \pi \subseteq 
Z(\scA)$ is nothing but a central extension, see Section \ref{sec4} 
for more on central extensions and centroids. \psp

\Proof (a) That $\pi_{\End}$ is an algebra homomorphism is 
well-known and easily seen. Since $\Ker \pi$ is an ideal, all left 
and right multiplication operators of $\scA$ leave $\Ker \pi$ 
invariant, whence  $\mult(\scA) \subseteq \End_\K(\scA;\Ker \pi)$. 
Also, for the left multiplication operator $L_x$ on $\scA$ we have 
$\pi\circ L_x = L_{\pi(x)} \circ \pi$ which shows $\pi_{\End}(L_x) = 
L_{\pi(x)}$. Since the analogous formula holds for the right 
multiplication operators and since $\pi_{\End}$ is an algebra 
epimorphism, we have (\ref{pimult}). Let $\chi \in \cent(\scA) \cap 
\End_\K(\scA;\Ker \pi)$. Then for $x,y \in \scA$, we have $ \bar 
\chi\big(\pi(x)\pi(y)\big) = \bar \chi \pi(xy) = \pi \chi (xy) = 
\pi\big( x\chi(y)\big) = \pi(x) (\bar \chi \pi(y)) = \pi(\chi(x)y)= 
\bar \chi(\pi(x))\pi(y)$,  which proves $\bar \chi \in \cent(\scB)$.

(b) If $\bar \chi = 0$ for $\chi \in \cent(\scA)\cap \hbox{\rm 
End}_\K(\scA;\Ker \pi))$, then $\chi(\scA) \subseteq \Ker \pi 
\subseteq \Ann_{\scA}(\scA)$.  So $\chi(xy) =\chi(x)y = 0$ for all 
$x,y \in \scA$, and because $\scA$ is perfect, it must be that $\chi 
= 0$.  

(c) It follows readily from $\pi \bigl(\Ann_{\scA}(\scA) \bigr) 
\subseteq  \Ann_{\scB}(\scB) = 0$ that $\Ker \pi=\Ann_{\scA}(\scA)$. 
By (a), $\pi_{\cent} : \cent(\scA) \to \cent(\scB)$ is then a 
well-defined algebra homomorphism, which is injective by (b). \qed 
\pse

 \subsection{Centroids of graded algebras}\label{secgrad}

We recall some concepts and results from the theory of graded 
algebras and graded modules (\cite[Sec.~11]{bou:A}). \psp 

 \begin{defn}\label{defgrad}  
 (1)  Let $\Lambda$ be an abelian group written additively.  An
 algebra $\scA$ over some base ring $\K$ is said to be
{\it $\Lambda$-graded} if $\scA=\bigoplus_{\lambda \in \Lambda}\, 
\scA^\lambda$ is a direct sum of $\K$-submodules $\scA^\lambda$ 
satisfying $\scA^\lambda \scA^\mu \subseteq \scA^{\lambda + \mu}$ 
for all $\lambda,\mu  \in \Lambda$. In this case, $\supp_{\scA} = 
\{\lambda \in \Lambda \mid \scA^\lambda \ne 0\}$ is called the  {\it 
support} of $\scA$, and the elements of $\scA^\lambda$ are said to 
be   {\it homogeneous of degree $\lambda$}.   A subalgebra (or 
ideal) $\scB$ of $\scA$ is  {\it graded} if $\scB=\bigoplus_{\lambda 
\in \Lambda}\, (\scB\cap \scA^\lambda)$.    Then $\scA$ is  {\it 
graded-simple} if $\scA^{(1)} \ne 0$,  and every graded ideal $\scB$ 
of $\scA$ is {\it trivial}, i.e., $\scB=0$ or $\scB=\scA$.  

(2) A $\Lambda$-graded unital associative algebra $\scA$ is said to 
be a {\it division-graded} algebra if every nonzero homogeneous 
element of $\scA$ is invertible. \end{defn} \psp

When $\scA$ is a division-graded associative algebra, $\supp_{\scA}$ 
is a subgroup of $\Lambda$; \ $\scA^0$ is division algebra;  and 
$\scA$ is a crossed product algebra \break $\scA = \allowbreak 
\scA^0 * \supp_\scA$.   Conversely, every crossed product algebra 
over a division algebra is a division-graded associative algebra.  
In particular, a commutative associative division-graded algebra 
$\scA$ is the same as a twisted group ring $\mathbb 
E^t[\supp_{\scA}]$ for $\mathbb E=\scA^0$ (see for example,  
\cite[Sec.~1]{P}).  \psp

Now let  $\scB$ be a $\Lambda$-graded unital associative 
$\K$-algebra.   A left $\scB$-module $\scM$ is {\it 
$\Lambda$-graded} if $\scM$ is a direct sum of $\K$-submodules,  
$\scM=\bigoplus_{\lambda \in \Lambda}\, \scM^\lambda$, such that 
$\scB^\lambda \scM^\mu \subseteq \scM^{\lambda + \mu}$  for all 
$\lambda,\mu \in \Lambda$.    In this case, we denote by 
$\End_\scB(\scM)^\lambda$  the $\K$-submodule  of all $f \in 
\End_\scB (\scM)$ satisfying $f\scM^\mu \subseteq \scM^{\lambda + 
\mu}$ for all $\mu \in \Lambda$, and we set
$$ \grEnd_\scB(\scM) = \bigoplus_{\lambda \in \Lambda} \End_\scB(\scM)^\lambda.
$$
This is a $\Lambda$-graded associative subalgebra of 
$\End_\scB(\scM)$ such that  $\scM$ is canonically a 
$\Lambda$-graded left module over $\grEnd_\scB(\scM)$.   In general 
$\grEnd_\scB(\scM)$ is a proper subalgebra of $\End_\scB(\scM)$. 
However, the following is proven in \cite[Sec.~11.6, Rem.]{bou:A}: 
\psp

\begin{lem}\label{remarque}
If $\scM$ is a finitely generated graded $\scB$-module,  then  
$\grEnd_\scB(\scM)$  $= \End_\scB(\scM)$.
\end{lem}
\psp

A $\Lambda$-graded $\scB$-module $\scM$ is {\em graded-irreducible} 
if the only graded $\scB$-sub-modules $\scN =\bigoplus_{\lambda\in 
\Lambda} \, (\scN \cap \scM^\lambda)$ are the trivial submodules 
$\scN= 0$ and $\scN=\scM$. A straightforward adaptation of the usual 
proof of Schur's Lemma gives the graded version below, in which the 
equality $\End_\scB(\scM) = \grEnd_\scB(\scM)$ follows from Lemma 
\ref{remarque}. \psp

\begin{lem}\label{Schur}
Let $\scB$ be a $\Lambda$-graded associative algebra and let $\scM$ 
be a $\Lambda$-graded $\scB$-module which is graded-irreducible. 
Then $\End_\scB(\scM) = \grEnd_\scB(\scM)$ is a division-graded 
algebra.
\end{lem} 
\psp

We now apply the above to a $\Lambda$-graded algebra $\scA$ over 
$\K$.   The multiplication algebra $\mult(\scA)$ is a graded 
subalgebra of the $\Lambda$-graded algebra $\grEnd_\K\scA$, and 
$\scA$ is a $\Lambda$-graded $\mult(\scA)$-module.     Then since 
$\cent(\scA) = \End_{\mult(\scA)} \scA$, we have 
\begin{eqnarray*}\label{} 
 \grcent(\scA)  &=& \grEnd_{\mult(\scA)} \scA 
   = \bigoplus_{\lambda \in \Lambda}\, \cent(\scA)^\lambda,  \quad\hbox{where}\\
 \cent(\scA)^\lambda &=& \cent(\scA) \cap \End_\K(\scA)^\lambda. 
\label{}
\end{eqnarray*}
By Lemma \ref{remarque} we see that 
\begin{equation}
\grcent(\scA) = \cent(\scA)\ \ \hbox{\it  if $\scA$ is a finitely 
generated $\mult(\scA)$-module.} \label{fingen}
\end{equation}
The graded version of Schur's Lemma now yields the following result. 
\psp

\begin{prop}\label{gradprop}
Let $\scA$ be a $\Lambda$-graded $\K$-algebra that is graded-simple. 
Then $\grcent(\scA)=\cent(\scA)$ is a division-graded commutative 
associative algebra, i.e., a twisted group ring $\mathbb 
E^t[\Gamma]$ over an extension field $\mathbb E = \cent(\scA)^0$ of 
$\K$  where $\Gamma=\{\lambda \in \Lambda  \mid \cent(\scA)^\lambda 
\ne 0\}$ is a subgroup of $\Lambda$.   Moreover, for every nonzero 
homogeneous $a\in \scA$, the evaluation map
$$\ev_a : \cent(\scA) \to \scA, \quad  \chi \mapsto \chi(a)
$$
is an injective $\mult(\scA)$-module map of degree $0$.
\end{prop}

\Proof  By assumption, $\scA$ is a graded-irreducible 
$\mult(\scA)$-module. {F}rom  Lemma \ref{Schur}, we   know that 
$\grcent(\scA)=\cent(\scA)$ is a division-graded associative 
algebra. It is commutative,  since $\scA^{(1)}$ is a graded ideal 
and hence $\scA$ is perfect.    It follows from the definitions that 
$\ev_a$ is a $\mult(\scA)$-module map of degree $0$.   It is 
injective since the $\mult(\scA)$-submodule generated by $a$ is all 
of $\scA$, i.e., $\mult(\scA). a = \scA$.    \qed \pse

\subsection{Centroids of tensor products and loop algebras}
 In the following all tensor products will be over a unital 
commutative associative ring $\K$. Let $\scA$ and $\scB$ be 
$\K$-algebras. There exists a unique $\K$-algebra structure on $\scA 
\ot \scB$ satisfying $(a_1 \ot b_1)(a_2 \ot b_2) = (a_1a_2) \ot 
(b_1b_2)$ for $a_i \in \scA$ and $b_i \in \scB$. Also, for $f\in 
\End_\K(\scA)$ and $g\in \End_\K(\scB)$ there exists a unique map $f 
\wot g\in \End_\K(\scA \ot \scB)$ such that $(f \wot g)(a\ot b) = 
f(a) \ot g(b)$ for all $a\in \scA$ and $b\in \scB$. The map $f\wot 
g$ should not be confused with the element $f\ot g$ of the tensor 
product $\End_\K(\scA) \ot \End_\K(\scB)$. Of course, one has a 
canonical map \begin{equation} \label{omega}
  \om : \End_\K (\scA) \ot \End_\K(\scB) \to  \End_\K(\scA \ot 
  \scB) :  \   f\ot g \mapsto f\wot g, 
 \end{equation}
It is straightforward to see that if $\chi_\scA \in \cent(\scA)$ and 
$\chi_\scB \in \cent(\scB)$  then $\chi_\scA \wot \chi_\scB \in 
\cent(\scA \ot \scB)$, and so $\cent(\scA ) \wot \cent(\scB) 
\subseteq \cent(\scA \ot \scB)$ where $\cent(\scA) \wot \cent(\scB)$ 
is the $\K$-span of all endomorphisms $\chi_\scA \wot \chi_\scB$. 

We will say that $\chi \in \cent(\scA \ot \scB)$ {\it has finite 
$\scA$-image} if for every $b\in \scB$,  there exist finitely many 
$b_1, \ldots , b_n \in \scB$ such that $\chi(\scA \ot \K b) 
\subseteq \scA\ot \K b_1 + \cdots + \scA \ot \K b_n$. It is easily 
seen that    
\begin{equation}\label{finim}
\cent(\scA) \wot \cent(\scB) \subseteq \{ \chi \in \cent(\scA \ot 
\scB) \mid \chi \hbox{ has finite $\scA$-image} \}. 
\end{equation}
We mention that for a unital $\K$-algebra $\scB$,  a $\chi \in 
\cent(\scA \ot \scB)$ has finite $\scA$-image as soon as $\chi(\scA 
\ot 1) \subseteq \scA\ot \K b_1 + \cdots + \scA \ot \K b_n$ for 
suitable $b_i\in \scB$.

\begin{prop}\label{Elemgeneral}
Let ${\scA}$ be a perfect $\K$-algebra and let ${\scB}$ be a unital 
$\K$-algebra that is free as a $\K$-module.   Then
\begin{itemize} 
 \item[{\rm (a)}] $\scA \ot \scB$ is perfect. 
 \item[{\rm (b)}] Every $\chi \in \cent(\scA \ot \scB)$ has finite 
 $\scA$-image if either one of the following conditions holds:
\begin{itemize}
  \item[{\rm (b.1)}] $\scA$ is finitely generated as a $\mult(\scA)$- or as a
  $\cent(\scA)$-module, or
   \item[{\rm (b.2)}] $\scA$ is central and a torsion-free $\K$-module. 
 \end{itemize}
  \item[{\rm (c)}] If\/ $\cent(\scA)$ is a free $\K$-module and the map 
    $\om$ of (\ref{omega}) is injective, then 
$$ \cent(\scA) \wot \cent(\scB) =  \{ \chi \in 
\cent(\scA \ot \scB) \mid \chi \hbox{ has finite $\scA$-image} \}. 
$$ 
\end{itemize}
\end{prop}

\Proof  (a) Since $\scA$ is perfect,  any $a\ot b \in \scA \ot \scB$ 
can be written as a finite sum $a\ot b = \sum_i (a'_i a_i'') \ot b = 
\sum_i (a'_i \ot b)(a_i'' \ot 1)$, where $1$ denotes the identity 
element of $\scB$. This implies (a).     

(b.1) Set $\scM = \mult(\scA)$, and observe that $\scM \ot \Id 
\subseteq \mult(\scA \ot \scB)$ since $\scB$ is unital. Suppose 
$\scA = \scM a_1 + \cdots + \scM a_n$ for $a_1, \ldots , a_n \in 
\scA$ and fix $\chi \in \cent(\scA \ot \scB)$ and $b\in \scB$. There 
exist finite families $(a_{ij}) \subseteq \scA$ and $(b_{ij}) 
\subseteq \scB$ such that  $\chi(a_i \ot b ) = \sum_j a_{ij} \ot 
b_{ij}$ for $1\le i \le n$, and hence $\chi(\scA \ot b) = \sum_{i,j} 
\, \chi\big( (\scM \wot \Id) (a_i \ot b) \big) = \sum_{i,j} (\scM 
\wot \Id)(a_{ij} \ot b_{ij}) \subseteq \sum_{i,j} \, \scA \ot 
b_{ij}$. 

By (a) the centroid $\cent(\scA \ot \scB)$ is commutative. The same 
argument as above with $\scM$ replaced by $\cent(\scA)$ then shows 
that every $\chi \in\cent(\scA \ot \scB)$ has finite $\scA$-image if 
$\scA$ is a finitely generated $\cent(\scA)$-module. 

(b.2) and (c): (This part of the proof is inspired by 
\cite[Lem.~1.2]{BM}.) Let $\{b_r\}_{r\in \RR}$ be a basis of $\scB$, 
and let $\chi \in \cent(\scA \ot \scB)$. We define $\chi_r \in 
\End_\K (\scA)$ by 
 \begin{equation} \label{chidef} \chi(a \ot 1) = \sum_{r\in \RR} \, 
\chi_r(a) \ot b_r\, . \end{equation} For $a_1, a_2 \in \scA$ we then 
have  
\begin{eqnarray*}
 \chi(a_1a_2\ot 1) &=& \sum_{r\in \RR} \, \chi_r (a_1a_2) \ot b_r \\
       &=& \chi\big( (a_1 \ot 1) (a_2 \ot 1)\big) 
            = \big(\chi(a_1 \ot b)\big) (a_2 \ot 1)\\
    &=& \sum_{r\in \RR} \, \chi_r(a_1) a_2 \ot b_r\\
       &=& (a_1 \ot 1) \chi(a_2\ot b) = \sum_{r\in \RR}
                     \,a_1\chi_r(a_2) \ot b_r, \\
 \end{eqnarray*}
whence $\chi_r(a_1a_2) = \chi_r(a_1)a_2 = a_1 \chi_r(a_2)$ for all 
$a_i \in \scA$, so all $\chi_r\in \cent(\scA)$. We can now finish 
the proof of (b.2): By assumption there exist scalars $x_r\in \K$ 
such that $\chi_r = x_r\Id$, hence $\chi(a\ot 1) = \sum_{r\in \RR} 
x_r a \ot b_r$. Fix $a\in \scA$. Then almost all $x_ra = 0$, so 
almost all $x_r=0$, which in turn implies that $\chi$ has finite 
$\scA$-image.   

We continue with the proof of (c). Because of (\ref{finim}),  we 
only need to prove the inclusion from right to left. So suppose that 
$\chi \in\cent(\scA \ot \scB)$ has finite $\scA$-image. Then there 
exists a finite subset $\FF\subseteq \RR$ such that  (\ref{chidef}) 
becomes
$$  \chi(a\ot 1) = \sum_{r\in \FF} \, \chi_r (a) \ot b_r\, .
$$
For $a_1, a_2 \in \scA$ and $b\in \scB$,  we then get 
$\chi(a_1a_2\ot b) = \chi\big( (a_1\ot 1)(a_2\ot b)\big) = 
\big(\chi(a_1 \ot 1)\big)(a_2 \ot b) = \sum_{r\in \FF} \, \chi_r 
(a_1)a_2 \ot b_r b = \sum_{r\in \FF} \, \chi_r (a_1a_2) \ot b_r b$. 
Since $\scA$ is perfect,  this implies 
\begin{equation}\label{abpform}
 \chi = \sum_{r\in \FF} \, \chi_r {\wot} \lam_r  
\end{equation}
where $\lam_r$ is the left multiplication in $\scB$ by $b_r$. 

Let $\{ \psi_s \mid s\in{\mathfrak S}\}$ be a $\K$-basis of 
$\cent(\scA)$. Then there exists a finite subset $\GG\subseteq 
{\mathfrak S}$ and scalars $x_{rs}\in \K$ $(r\in \FF, s\in \GG)$ 
such that $\chi_r = \sum_{s\in \GG} \, x_{rs}\psi_s$. We then get 
$\chi = \sum_{r\in \FF} \sum_{s\in \GG} \, x_{rs} \psi_s \wot \lam_s 
= \sum_{s\in \GG}\, \psi_s \wot \phi_s$ for  $\phi_s = \sum_{r\in 
\FF} x_{rs}\lam_r$, and it remains to show that $\phi_s \in 
\cent(\scB)$. For $a_i \in \scA$, $b_i \in \scB (i=1,2)$ we have 
\begin{eqnarray*}
\sum_{s\in \GG}  \, \psi_s(a_1a_2) \ot \phi_s (b_1b_2) 
  &=& \chi (a_1a_2 \ot b_1b_2) = \big(\chi(a_1\ot b_1)\big)(a_2\ot b_2) \\
  &=& \sum_{s \in \GG}  \psi_s(a_1)a_2\ot \phi_s(b_1)b_2 \\ 
  &=&\sum_{s\in \mathfrak \GG} \, \psi_s(a_1a_2) \ot \phi_s(b_1)b_2\,.
\end{eqnarray*}
Because $\scA^{(1)}=\scA$, this implies $\sum_{s\in \GG} \, 
\psi_s(a) \ot (\phi_s(b_1b) - \phi_s(b_1)b)=0$ for all $a\in \scA$ 
and $b, b_1 \in \scB$, and then $\sum_{s\in \GG}\, \psi_s \wot \mu_s 
= 0$ where $\mu_s \in \End_\K(\scB)$ is defined as $\mu_s(b) = 
\phi_s(b_1b) - \phi_s(b_1)b$. Since by assumption $\om$ is 
injective, we also have $\sum_{s\in \GG}\, \psi_s \ot \mu_s = 0$. So 
by the linear independence of the $\psi_s$, we see that $\mu_s = 0$.  
But then $\psi_s \in \cent(\scB)$ follows, and hence $\chi \in 
\cent(\scA) \wot \cent(\scB)$.   \qed 

\psp 

\begin{rem} {\rm Let $\scA$ and $\scB$ be $\K$-algebras such that 
$\scA$ is a finitely generated $\mult(\scA)$-module and $\scB$ is  
unital commutative associative $\K$-algebra which is free as a 
$\K$-module.  It is shown in \cite[Lemma~2.22]{ABP3} that then 
$$\cent(\scA) \ot \cent(\scB)\, {\buildrel \cong \over \longrightarrow}\, 
   \cent(\scA) \wot \cent(\scB) = \cent(\scA \ot \scB).$$
Indeed,  equation (\ref{chidef}), (which in fact holds even if 
$\scA$ is not perfect), together with (b.1) of Proposition 
\ref{Elemgeneral} shows that $\cent(\scA) \wot \cent(\scB) = 
\cent(\scA \ot \scB)$. The isomorphism $\cent(\scA) \ot 
\cent(\scB)\cong \cent(\scA) \wot \cent(\scB)$ follows from the 
linear independence of the set $(b_r)_{r\in \RR}\,$. }
\end{rem} \psp
 
{F}rom now on, unless explicitly stated otherwise,  all algebras 
will be over some field $\F$. For easier reference we state the 
following consequence of Proposition \ref{Elemgeneral}. 

\begin{cor}\label{Elem}
 Let ${\scA}$ and ${\scB}$ be algebras over a field $\F$ such
 that ${\scA}$ is perfect and ${\scB}$ is unital.
Then  
 \begin{eqnarray*}
 \cent(\scA) \ot \cent(\scB)  &\cong& \cent(\scA) \wot \cent(\scB) \\
 &=& \{ \chi \in  \cent(\scA \ot \scB) \mid \chi \hbox{ has finite $\scA$-image} 
   \}.
 \end{eqnarray*}
Moreover, we have 
 $$ \cent(\scA) \wot \cent(\scB) = \cent (\scA \ot \scB)$$
in either one of the following cases:
\begin{itemize}
 \item[{\rm (a)}] $\scA$ is finitely generated as $\mult(\scA)$-module, 
    e.g., $\dim_{\F}{\scA} < \infty$, or
 \item[{\rm (b)}] $\cent({\scA}) = \F \Id$, in which
 case $\cent({\scA} \ot {\scB}) = \Id \ot \cent({\mathcal B})$. 
\end{itemize}
\end{cor}

\Proof Over fields the map $\om$ is injective (\cite[Sec. 7.7, 
Prop.~16]{bou:A}). All assertions then follow from Proposition 
\ref{Elemgeneral}.   

 \psp

\begin{rem} \label{remloo} {\rm When $\g$ is a simple Lie algebra with $\cent(\g) = \F \Id$ (which
is always true when $\F$ is algebraically closed and $\g$ is 
finite-dimensional) Melville \cite[3.2, 3.4, 3.5]{Me}  has shown 
$\cent(\g \ot \mathcal B) = \F \Id \ot \cent(\mathcal B)$ for 
$\mathcal B = \F[x_1,\dots,x_n]$ a polynomial ring or an ideal 
$\mathcal B=t^m\F[t^n]$ of the polynomial ring $\F[t]$. Actually the 
proof given in \cite[3.2]{Me}  works for any unital commutative 
associative  $\F$-algebra $\mathcal B$.  Part (b) of 
Corollary~\ref{Elem} was proven by Allison-Berman-Pianzola 
(\cite[Lem.~4.2]{ABP2}) under the assumption that  $\mathcal B$ is a 
unital commutative associative $\F$-algebra.}  \end{rem} \psp

\begin{exam} {\rm   Suppose ${\scA}$ is the centreless Virasoro
Lie algebra (often called the Witt algebra).  Thus, ${\scA}$ has a 
basis consisting of the elements $\{a_i \mid i \in \ZZ\}$ and 
multiplication given by $[a_i, a_j] = (j-i)a_{i+j}$.   Suppose $\chi 
\in \cent(\scA)$.  Then $j\chi(a_j) = \chi\bigl([a_0,a_j]\bigr) =  
[a_0,\chi(a_j)]$.   Since $\{x \in {\scA} \mid [a_0,x] = j x\} = \F 
a_j$, we see that $\chi(a_j) = \lambda_j a_j$ for some scalar 
$\lambda_j$. But then $\lambda_{i+j} (j-i)a_{i+j} = 
\chi\bigl([a_i,a_j]\bigr)= [a_i,\chi(a_j)] = (j-i)\lambda_ja_{i+j}$ 
(with $j=0$)  shows that $\lambda_i = \lambda_0$ for all $i$.  Hence 
$\cent({\scA}) = \F \Id$ --- this can also be derived easily from 
Proposition \ref{toral} below.   Applying part (b) of the previous 
proposition, we obtain
 that  $\cent({\scA}
\ot {\scB}) =  \Id \ot \cent({\scB})$ for any unital algebra 
${\scB}$. }
\end{exam} \psp

\begin{cor}\label{Elemcor}
Let $\scA$ be a central, perfect $\F$-algebra,  $\scB$ be a unital 
commutative associative $\F$-algebra,  and $\scC$ be a unital 
subalgebra of $\scB$ such that $\scB$ is a  free $\scC$-module with 
a $\scC$-basis containing the identity element $1$ of $\scB$.  
Suppose further that $\scL$  is $\scC$-form of $\scA\ot_\F \scB$, 
i.e., a subalgebra of the $\scC$-algebra $\scA\ot_\F\scB$ such that 
$\scL \ot_\scC \scB \cong \scA \ot_\F \scB$. Then $\cent_\F(\scL) = 
\scC\, \Id$.
\end{cor}

\Proof  It follows from the assumptions that $\scC\, \Id  \subseteq  
\cent_\F(\scL)$, where $a \otimes b \mapsto a \otimes cb$ for all $c 
\in \scC$.  Conversely, let $\chi \in \cent_\F(\scL)$ and extend 
$\chi$ to a $\F$-linear map $\widetilde \chi$ on $\scL \ot_\scC\scB$ 
by setting  $\widetilde \chi = \chi \ot \Id$. Since $\scL \ot_\scC 
\scB \cong \scA \ot_\F \scB$ and $\cent_\F(\scA \ot_\F\scB) \cong  
\Id \ot_\F \scB$ by Corollary~\ref{Elem},  it follows that there 
exists $c\in \scB$ such that  $\widetilde \chi $ is given by $x 
\ot_\scC b \mapsto x\ot_\scC cb$ for all $x\in \scL$, $b \in \scB$. 
Let $(b_i)_{i\in \scI}$ be a $\scC$-basis of $\scB$ containing $1$, 
say $b_0 = 1$.     We can write $c$ in the form $c=\sum_i c_i b_i$ 
for unique $c_i\in C$. Then $\widetilde \chi(x\ot_\scC 1) = 
x\ot_\scC c = \sum_i x\ot_\scC c_ib_i = \sum_i xc_i \ot_\scC b_i$. 
Since $\widetilde \chi(x\ot_\scC 1)= \chi(x) \ot_\scC 1 \in 
\scL\ot_\scC 1$,  it follows that  $xc_i = 0$ for $i\ne 0$ and 
$\chi(x) = c_0x$, i.e., $\chi = c_0\,\Id$.   Thus, $\cent_\F(\scL) 
\subseteq \scC\, \Id$, and we have the desired conclusion.   \qed 
\psp

\begin{rem}\label{laurent}  {\rm The assumptions on $\scB$ and $\scC$ are fulfilled
for example when $\scB =\F[t,t^{-1}]$, a Laurent polynomial ring, 
and $\scC=\F[t^m, t^{-m}]$ for some positive $m\in \NN$.    Indeed, 
in this case $\{t^i \mid  0\le i < m\}$ is a $\scC$-basis of $\scB$ 
as required.      Let $\scA$ be a Lie algebra and $\sigma$ be an 
automorphism of $\scA$ of order $m$.   Assume $\zeta \in \F$ is a 
primitive $m$th root of unity, and form the loop algebra 
$L(\scA,\sigma) = \bigoplus_{i\in \ZZ} \scA_i \ot_\F \F t^i$  where 
$\scA_i$ is the $\zeta^i$-eigenspace of $\sigma$. That 
$L(\scA,\sigma)$ is indeed a $\scC$-form of $\scA \ot_\F \scB$ is 
shown by  Allison-Berman-Pianzola \cite[Thm.~3.6]{ABP2}. In this 
particular situation,  the result  $\cent(L(\scA,\sigma)) = \scC\, 
\Id$ from  Corollary \ref{Elemcor}
 can be found in  \cite[Lem.~4.3 (d)]{ABP2}.   
When  $\scA$ is taken to be a finite-dimensional split simple Lie 
algebra over a field $\F$ of characteristic $0$, the loop algebra 
$L(\scA,\sigma)$  is an example of a centreless core of an extended 
affine Lie algebra.   Later in Corollary \ref{ealacent},  we will 
see that the cores of extended affine Lie algebras are always 
central.}
\end{rem}  \psp

\begin{cor}\label{autocent}
Let $\scA$ be an algebra over a ring $\K$ and set $\scC = 
\cent_\K(\scA)$.
\begin{itemize}
\item[{\rm (a)}]  Every $\K$-linear automorphism $f$ of $\scA$ determines
a $\K$-linear automorphism  $ f_{\scC} :  \scC \to \scC, \quad  \chi 
\mapsto f_{\scC}(\chi) =f\circ \chi\circ f^{-1}$ of the associative 
$\K$-algebra $\scC$.  The map
$$
  \ga : \Aut_\K(\scA) \to \Aut_\K(\scC), \quad f\mapsto f_{\scC}$$ is a group homomorphism
whose kernel is $\Aut_\scC(\scA)$, the $\scC$-linear automorphisms 
of $\scA$.

\item[{\rm (b)}]  Let $\scA$  be a perfect, central algebra over a field $\F$ and  $\scB$ be a unital 
commutative associative $\F$-algebra.  Then, after identifying 
$\Aut_\F(\scB)$ with a subgroup of $\Aut_\F(\scA\ot_\F \scB)$ via 
$g\mapsto \Id \ot  g$, we have
$$   \Aut_\F(\scA \ot_\F \scB) = \Aut_\scB (\scA \ot_\F \scB)
\rtimes \Aut_\F(\scB) \quad (\hbox{semidirect product}).
$$  \end{itemize}
\end{cor} 

\Proof (a) is straightforward. For (b) we note that $\cent(\scA\ot 
\scB) = \Id \ot \scB \cong \scB$ by Corollary~\ref{Elem}.  Every 
$g\in \Aut_\F(\scB)$ extends to an automorphism $\xi(g) = \Id \ot g$ 
of $\scA \ot \scB$. The map $\xi : \Aut_\F(\scB) \to \Aut_\F(\scA)$ 
is a group homomorphism which satisfies $(\ga \circ \xi)(g) = g$, 
i.e., $\xi$ is a section of $\gamma$. The claim then follows from 
standard facts in group theory. \qed \psp

\begin{rems}\label{autorem} {\rm For special choices of $\scA$ and $\scB$, the automorphism
group $\Aut_\F(\scA \ot_\F \scB)$ has been investigated by several 
authors;  for example, by Benkart-Osborn \cite{BO}  when $\scA$ is 
the algebra of $n \times n$-matrices over $\F$ and $\scB$ is an 
arbitrary algebra with an Artinian nucleus, or by Pianzola \cite{Pi} 
when $\scA$ is a finite-dimensional split simple Lie algebra over a 
field $\F$ of characteristic $0$ and $\scB$ is an integral domain 
with trivial Picard group and with a maximal ideal $\mathfrak m$ 
which satisfies $\scB / {\mathfrak m} \cong \F$.  }
\end{rems}  \psp

  An analogous result holds for derivations;  we leave the proof
as an exercise to the reader.

\begin{cor}\label{dercent}
Let $\scA$ be an algebra over a ring $\K$ and set $\scC = 
\cent_\K(\scA)$.  
\begin{itemize} 
\item[{\rm (a)}]  Every $\K$-linear derivation $d$ of $\scA$  determines 
a $\K$-linear derivation $ d_{\scC} : \scC \to \scC, \quad \chi 
\mapsto d_{\scC} (\chi) := [d, \chi] = d\circ\chi -\chi\circ d$ of 
 $\scC$.   The map
$$ \delta : \Der_\K(\scA) \to \Der_\K(\scC), \quad
         d\mapsto d_{\scC}
$$ is a  $\K$-linear Lie algebra homomorphism
whose kernel is $\Der_\scC(\scA)$, the $\scC$-linear derivations of 
$\scA$. 

\item[{\rm (b)}] Let $\scA$  be a perfect, central algebra over a field $\F$ and  $\scB$ be a unital 
commutative associative $\F$-algebra.  Then, after identifying  
$\Der_\F(\scB)$ with a subalgebra of $\Der_\F(\scA\ot_\F \scB)$ via 
$e\mapsto \Id \ot e$, we have
$$   \Der_\F(\scA \ot_\F \scB) = \Der_\scB (\scA \ot_\F\scB)
\rtimes \Der_\F(\scB) \quad (\hbox{semidirect product}).
$$\end{itemize}
\end{cor}
\ms

 \begin{rems} \label{BenMoo} {\rm The result in
Corollary \ref{dercent}(b) complements \cite[Thm.~1]{BM},  which 
describes $\Der_\F(\scA \ot_\F \scB)$   when $\scA$ is a 
finite-dimensional perfect $\F$-algebra and $\scB$ is as above: 
\begin{eqnarray}
\Der_\F(\scA \ot_\F \scB) = \big(\Der_\F(\scA) \ot_\F \scB \big) \; 
\oplus \; \big(\cent_\F(\scA)\ot_\F \Der_\F(\scB)\big). \label{not}
\end{eqnarray}
We note that (\ref{not}) is not a semidirect product in general. 
}\end{rems} \pse

\section{Centroids of Lie algebras with toral subalgebras}\label{toralsec}

\subsection{A general result}\label{toralsecgen} In this section,  $\scL$ is a Lie algebra over
some field $\F$, which will be assumed of characteristic $0$ from 
Section \ref{Kac-Moody} on.    Recall that a subalgebra $\h$  is  a 
{\it toral subalgebra} of $\scL$  if $\scL= \bigoplus_{\al \in 
\h^*}\,\scL_\al,$ where $\scL_\al = \{ x\in \scL \mid [h,x]=\al(h)x 
\hbox{ for all } h\in \h \}$. Necessarily  $\h$ is abelian, and $ 
[\scL_\al \,,\, \scL_\beta] \subseteq \scL_{\al + \beta}$ for 
$\al,\beta\in \h^*$. If $\scL_\al \neq 0$, then $\alpha$ is a {\it 
weight} (relative to $\h$) and $\supp_\scL=\{\al\in \h^* \mid 
\scL_\al \ne 0\}$  is the set of weights.  
\begin{prop}\label{toral}
Let $\scL$ be a Lie algebra with a toral subalgebra $\h$.
\begin{itemize} \item[{\rm (a)}] If $\chi \in \cent(\scL)$, then 
$ \chi(\scL_\al) \subseteq \scL_\al
    \hbox{ for all } \al \in \h^*$,\ 
 $ \chi(\h) \subseteq  Z(\scL_0)$, \ and 
 $\chi$ is uniquely determined by its restriction to $\h$. 

\item[{\rm (b)}]  Let $\scJ$ be a $\cent(\scL)$-invariant ideal,  and
suppose there exists $0\ne \al\in \h^*$ such that $\dim (\scJ\cap 
\scL_\al) = 1$ and the ideal of $\scL$ generated
 by $\scJ\cap \scL_\al$ is $\scJ$. Then
\begin{eqnarray}
 \cent(\scL) &=& \F\Id \,\oplus\, \{\psi\in \cent(\scL)
              \mid \psi(\scJ)=0\} \nonumber \\
   &\cong& \F\Id  \, \oplus \, \Hom_{\scL/\scJ} (\scL/\scJ, \scC_\scL(\scJ))
   \label{toral3}
 \end{eqnarray}
where $\scC_\scL(\scJ)$ is the centralizer of $\scJ$ in $\scL$. In 
particular, if $\scJ=\scL^{(1)}$ satisfies the assumptions above, 
then
 \begin{eqnarray}
 \cent(\scL) &=& \F\Id  \, \oplus\, \{ \psi \in \End_\F(\scL)
               \mid  [\psi(\scL),\scL]=0=\psi(\scL^{(1)})\}   \label{toral4}\\
               &=& \F\Id  \, \oplus\, \scV(\scL^{(1)}) \nonumber
 \\
  &\cong& \F\Id  \,\oplus\, \Hom_\F\big(\scL/\scL^{(1)},  \scC_\scL(\scL^{(1)})\big). \nonumber
 \end{eqnarray}  \end{itemize} 
   \end{prop}
 
\Proof  (a) For all $h\in \h$ and $x_\al \in \scL_\al$ we have 
$\al(h) \chi(x_\al) = \chi[h,x_\al] = [h,\chi(x_\al)]$ which implies 
$\chi(\scL_\al) \subseteq \scL_\al$. Also $[\chi(\h), 
\scL_0]=\chi[\h,\scL_0]=0$, proving $\chi(\h)$ is central in 
$\scL_0$.  For $0\ne \al$ there exists $t_\al \in \h$ such that 
$\al(t_\al) = 1$. The last claim then follows from 
$\chi(x_\al)=\chi[t_\al,x_\al] = [\chi(t_\al),x_\al]$.

(b) Suppose $\chi\in \cent(\scL)$.   Our assumptions imply that 
there exists a scalar $\xi$ such that $\chi\, |_{\scJ\cap\scL_\al} = 
\xi\Id$. Thus $\scJ \cap \scL_\al$ is contained in the kernel of 
$\psi=\chi - \xi\Id$, which is an ideal of $\scL$  since $\psi \in 
\cent(\scL)$.  As $\scJ\cap\scL_\al$ generates $\scJ$,  we have 
$\psi\, |_\scJ =0$. This proves $\cent(\scL)$ is contained in the 
right-hand side. The other direction is obvious. The second part of 
(\ref{toral3}) and the statement concerning $\scL^{(1)}$ follow from  
Lemma \ref{easy}.
 \qed

\begin{cor}\label{toralcor} Let $\scL$ be a Lie algebra with a toral subalgebra
$\h$ and suppose that $\scL^{(1)}$ is generated by elements $e_i, 
f_i, (1\le i\le n)$ such that the following conditions hold:

\begin{enumerate}

\item[\rm (i)] $\al_i^\vee =[e_i,f_i]\in \h$, and these elements
act on the generators $e_j, f_j$ by
$$ [\al_i^\vee, e_j]=a_{i,j}e_j \quad\hbox{and} \quad
[\al_i^\vee, f_j]=-a_{i,j}f_j \quad (1\le i,j\le n)
$$
where $a_{i,j} \in \F$.

\item[\rm (ii)] The matrix $\mathfrak A=(a_{i,j})$ is indecomposable in the
sense that after possibly renumbering the indices we have 
$a_{1,2}a_{2,3} \cdots a_{n-1,n} \ne 0$. Moreover, $a_{i,i} \ne 0$ 
for all $i=1,\dots, n$.

\item[\rm (iii)] For some $i$ and some $\al \in \h^*$,  we
have $\scL_\al=\F e_i$.
\end{enumerate}

\noindent Then $\cent(\scL) \cong \F \Id \oplus 
\Hom_\F\big(\scL/\scL^{(1)}, \scC_{\scL}(\scL^{(1)}) \big)$. In 
particular, $\scL$ is central if $\scL =\scL^{(1)}$ or 
$\scC_{\scL}(\scL^{(1)}) = 0$. 
\end{cor}

\Proof Fix $i$ as in (iii). By Proposition \ref{toral}\,(b),  it 
suffices to show that the ideal $\scJ$ of $\scL$ generated by $e_i$ 
contains all the generators $e_j, f_j$ of $\scL^{(1)}$. This follows 
from the indecomposibility of the matrix $\mathfrak A$  by upward 
and downward induction starting from $i$. Indeed, assume $e_j\in 
\scJ$. Relation (i) implies that $\al_j^\vee \in \scJ$ and then, 
since $a_{j,j}\ne 0$, that $f_j\in \scJ$. Moreover, $e_{j+1} = 
a_{j,j+1}^{-1}[\al_j^\vee,e_{j+1}] \in \scJ$ and similarly $f_{j+1} 
\in \scJ$. \qed \pse

\subsection{Centroids of Kac-Moody algebras}\label{Kac-Moody}
 
{From now on} we assume $\F$ is a field of characteristic $0$. 
Corollary~\ref{toralcor} applies to contragredient Lie algebras over 
$\F$ in the sense of \cite[Ch.~4]{MP}. In particular, it applies to 
Kac-Moody algebras.  We will elaborate on this important special 
case.  \psp

Assume ${\mathfrak A}: = (a_{i,j})_{i,j=1}^n$ is a generalized 
Cartan matrix of rank $\ell$.   Thus, $a_{i,i} = 2$ for 
$i=1,\dots,n$; \ the entries $a_{i,j}$ for $i \neq j$ are 
nonpositive integers;  and $a_{i,j} = 0$ if and only if $a_{j,i} = 
0$.    We will assume that the matrix $\mathfrak A$ is 
indecomposable as in Corollary \ref{toralcor}(ii).   A realization 
of $\mathfrak A$ is a triple $({\mathfrak h}, \Pi, \Pi^\vee)$ 
consisting of an $\F$-vector space $\mathfrak h$ and linearly 
independent subsets $\Pi = \{\alpha_1,\dots,\alpha_n\} \subset 
\mathfrak h^*$ and $\Pi^\vee = \{\alpha_1^\vee, \dots, 
\alpha_n^\vee\}\subset \mathfrak h$ such that $\langle \alpha_i 
\,|\, \alpha_j^\vee\rangle = a_{j,i}$ ($1 \leq i,j \leq n$), and 
$\dim {\mathfrak h} = 2n-\ell $.  The {\it Kac-Moody Lie algebra} 
${\tt g} : = {\tt g}(\mathfrak A)$ associated to $\mathfrak A$ is 
the Lie algebra over $\F$ with generators $e_i, f_i, (1 \leq i \leq 
n)$ and $\mathfrak h$, which satisfy the defining relations,

\begin{itemize}
\item [{\rm (a)}]   $[e_i,f_j] = \delta_{i,j}\alpha_i^\vee$
\qquad $1 \leq i,j \leq n$; \item [{\rm (b)}]   $[h,h'] = 0$ \qquad 
$h,h' \in \mathfrak h$; \item [{\rm (c)}]   $[h,e_i] = \langle 
\alpha_i\,|\, h\rangle e_i$ and $[h,f_i] = - \langle \alpha_i \,|\, 
h\rangle f_i$ \qquad $h \in \mathfrak h$, \ $i=1,\dots,n$; \item 
[{\rm (d)}]   $(\ad e_i)^{1-a_{i,j}}e_j = 0 = (\ad 
f_i)^{1-a_{i,j}}f_j$  \qquad $1 \leq i\neq j \leq n$.

\end{itemize}
 
The Lie algebra ${\tt g}$ is graded by the root lattice $Q: = 
\bigoplus_{i=1}^n \mathbb Z \alpha_i$, so that  ${\tt g} = 
\bigoplus_{\alpha \in Q}{\tt g}_\alpha$, where ${\tt g}_{\alpha} = 
\{x \in{\tt g} \mid [h,x] = \langle \alpha \,|\, h\rangle x$ for all 
$h \in \mathfrak h\}$ are the root (weight) spaces of the toral 
subalgebra $\h$. Condition (i) of Corollary \ref{toralcor} follows 
from  relations (a) and (c) above,  while (ii) holds since 
$a_{i,i}=2$ for all $i$ and since the Cartan matrix $\mathfrak A$ is 
assumed to be indecomposable. Finally, (iii) is well-known, see e.g. 
\cite[(1.3.3)]{K}.   Therefore, we obtain \psp

\begin{cor}\label{centKM}  Let  ${\tt g} = {\tt g}(\mathfrak A)$ be the  Kac-Moody
Lie algebra corresponding to the indecomposable generalized Cartan 
matrix $\mathfrak A$.   Then ${\tt g}^{(1)}$ is central and 
$\cent({\tt g}) = \mathbb F \Id  \oplus\Hom_{\mathbb F} \bigl({\tt 
g}/{\tt g}^{(1)}, \scC_{{\tt g}}({\tt g}^{(1)}) \bigr).$
\end{cor}  \psp

\begin{rems} \label{remKM} {\rm When $\mathfrak A$ is invertible, $\dim \mathfrak h = n$,  \ $\mathfrak h  = \bigoplus_{i=1}^n \mathbb F \alpha_i^\vee$,  and 
${\tt g} = {\tt g}^{(1)}$, so that  $\cent({\tt g}) = \mathbb F \Id$ 
in that case.   

When $\mathfrak A$ is an affine Cartan matrix (associated to an 
affine Dynkin diagram),  then $\mathfrak A$ has rank $n-1$.     The 
center is one-dimensional, spanned by $c$ say.  We may  suppose 
$\mathfrak h = \mathfrak h' \oplus \mathbb C d$ in this case, where 
$\mathfrak h' = \bigoplus_{i=1}^n \mathbb F \alpha_i^\vee$ and ${\tt 
g} = {\tt g}^{(1)} \oplus \F d$.  (Readers familiar with affine 
algebras will recognize $d$ as the degree derivation.)  Corollary~ 
\ref{centKM} then says that ${\tt g}^{(1)}$ is central and 
$\cent({\tt g}) = \mathbb F \Id  \oplus\Hom_{\mathbb F}\bigl(\mathbb 
F d, \mathbb F c\bigr)$ for each affine Kac-Moody Lie algebra.  }
\end{rems} \pse
 
\subsection{Centroids of Lie tori} 

First  we introduce the notion of a root graded Lie algebra.   We 
prove a  result concerning the centroid of a special class of root 
graded Lie algebras  but postpone giving the precise description of 
the centroid for general root graded Lie algebras  until  Section 
\ref{rootgrad}.     We then specialize considerations  to certain 
root graded Lie algebras called  Lie tori and describe their 
centroids.   Our rationale for doing this is that Lie tori play a 
critical role in the theory of extended affine Lie algebras (EALAs) 
-- they are precisely the cores of EALAs.     The centroid of the 
core is another  key ingredient in the structure of the EALA.  \psp

Let $\g$ be a finite-dimensional split  simple Lie algebra over a 
field $\F$ of characteristic 0 with root space decomposition $\g = 
\h \oplus \bigoplus_{\alpha \in \Delta_\g} \g_\alpha$ relative to a 
split Cartan subalgebra $\h$.  Such a Lie algebra is the 
$\F$-analogue of a finite-dimensional complex Lie algebra, and the 
set of roots $\Delta_\g$ is the  finite reduced root system.  Every 
finite irreducible root system $\Delta$  is one of the reduced root 
systems  $\Delta_\g$  or is a nonreduced root system BC$_r$  (see 
for example, \cite[Sec.~1.1]{bou:lieVI}).  \psp

\begin{defn}\label{rootgradE} Let $\scL$ be a Lie algebra over a field $\F$ of
characteristic $0$, and let $\De$ be a finite irreducible root 
system.  Then $\scL$ is said to be  {\em graded by the root system  
$\Delta$}  or  to be {\it $\De$-graded}  if
\begin{itemize}
 \item[{\rm (i)}]  $\scL$ contains as a subalgebra  a
finite-dimensional split ``simple'' Lie algebra $\g$, called the 
{\it grading subalgebra},  with split Cartan subalgebra $\h$; 

 \item[{\rm (ii)}] $\h$ is a toral subalgebra of $\scL$, and
 the weights of $\scL$ relative to $\h$ are in $\Delta \cup \{0\}$,  \begin{equation*} 
 \scL = \bigoplus_{\al \in \Delta \cup \{0\}}
 \scL_\al\end{equation*}
 \item[{\rm (iii)}]  $\scL_0 = \sum_{\al \in \Delta}\, [\scL_\al, \scL_{-\al}]$;
 \item[\rm (iv)] either $\De$ is reduced and equals the root system
$\De_\g$ of $(\g,\h)$  or $\De={\rm BC}_r$ and $\De_\g$ is of type  
B$_r$,  C$_r$,  or D$_r$.
\end{itemize}\end{defn}

The word simple is in quotes, because in all instances except two, 
$\g$ is a simple Lie algebra.  The sole exceptions are when $\Delta$ 
is of type BC$_2$,   $\Delta_\g$ is of type D$_2$ = A$_1 \times$ 
A$_1$,  and $\g$ is the direct sum of two copies of 
$\mathfrak{sl}_2$;  or when $\Delta$ is of type BC$_1$, $\Delta_\g$ 
is of type D$_1$, and $\g = \h$ is one-dimensional.

\psp The definition above is due to Berman-Moody \cite{BM} for the 
case $\De=\De_\g$.  The extension to the nonreduced root systems was 
developed by  Allison-Benkart-Gao in  \cite{ABG2}.    The 
$\De$-graded  Lie algebras for $\Delta$ reduced  have been 
classified up to central extensions by Tits \cite{T} for $\Delta =$ 
A$_1$ (see also \cite{BZ});  by Berman-Moody \cite{BeM} for $\Delta 
=$ A$_r$, ($r \geq 2$), D$_r$,  E$_6$, E$_7$, and E$_8$; \; by 
Benkart-Zelmanov  \cite{BZ} for $\Delta =$B$_r$, C$_r$, F$_4$ and 
G$_2$ (see also\cite{ABG1}  for C$_r$);\; and by Neher \cite{N1}, 
who studied Lie algebras  3-graded by a locally finite root system 
$\De$,  which in our setting means $\Delta \neq$ E$_8$, F$_4$, 
G$_2$, or BC$_r$.  Central extensions of Lie algebras graded by 
reduced root systems have been described by Allison-Benkart-Gao in 
\cite{ABG1}.    The Lie algebras graded by the root systems BC$_r$ 
have been classified in \cite{ABG2} for $r \geq 2$ and in \cite{BS} 
for $r=1$.    As a result, the $\Delta$-graded Lie algebras are 
determined completely up to isomorphism.  
 \psp
 
Examples of Lie algebras graded by a (not necessarily reduced) 
$\Delta$  include the affine Kac-Moody Lie algebras, the toroidal 
Lie algebras,  the intersection matrix Lie algebras introduced by 
Slodowy in his study of singularities, and the cores of EALAs, to 
name just a few (see  \cite[Exs. 1.16-1.23]{ABG2} for further 
discussion of examples).   Any finite-dimensional simple Lie algebra 
over a field $\mathbb F$ of characteristic 0 which has an 
ad-nilpotent element (or equivalently, by the Jacobson-Morozov 
theorem has  a copy of $\mathfrak {sl}_2$) is graded by a finite 
root system (see \cite{Se}).   Thus, the notion encompasses a 
diverse array of  important Lie algebras. \psp

\begin{defn}\label{defrc}
Let $\scL$ be a $\De$-graded Lie algebra with grading subalgebra 
$\g$, and let $\La$ be an abelian group.  We say that $\scL$ is {\it 
$(\De,\La)$-graded} if $\scL$ has a $\La$-grading
\begin{equation} \scL=\bigoplus_{\lam\in \La}\,
\scL^\lam  \quad\hbox{with}\quad  \g \subseteq \scL^0. \label{delam}
\end{equation} \end{defn} 
\psp

This notion was introduced and studied by Yoshii in \cite{Y2} where 
it was termed a {\it refined root grading of type $(\De,\La)$}.  In 
any $(\De,\La)$-graded Lie algebra,  each space $\scL^\lam$ is an 
$\h$-submodule, so $\scL^\lam = \bigoplus_{\al \in \De \cup \{0\}} 
\, \scL^\lam_\al$ for $ \scL^\lam_\al = \scL_\al \cap \scL^\lam$.   
Thus, $\scL$ has a grading by $\La \oplus \scQ(\De)$, where 
$\scQ(\De)$ is the root lattice of $\De$, 
\begin{equation}\scL =
\bigoplus_{\al\in\De \cup \{0\},\, \lam\in \La}\, \scL_\al^\lam 
\quad\hbox{with}\quad [\scL^\lam_\al\,,\, \scL^\mu_\beta] \subseteq 
\scL^{\lam + \mu}_{\al + \beta} \label{DeLamGrad}
\end{equation}  for $\al,\beta \in \De \cup \{0\}$ and $\lam,\mu \in \La$.
 Since the centre of any graded Lie algebra is a graded subspace,  
 $Z(\scL) = \bigoplus_{\lam \in \La}\, Z(\scL)^\lam \subseteq \bigoplus_{\lam \in \La}\, \scL^\lam_0$.  \psp
 
\begin{lem} \label{rootgradcen} {\rm (1)} Let $\scL$ be a $(\De,\La)$-graded Lie algebra. Then
$$ \cent(\scL)=  \bigoplus_{\lam\in \La}\,
\cent(\scL)^\lam
$$
is a $\La$-graded commutative associative algebra, where 
$\cent(\scL)^\lam$ is the subspace of centroidal transformations 
that are homogeneous of degree $\lam$ with respect to the 
$\La$-grading {\rm (\ref{delam})}. \smallbreak

{\rm (2)} If  the $\La$-graded Lie algebra $\scL$ is graded-simple, 
its centroid  $\cent(\scL)$ is a commutative associative 
division-graded algebra,  hence it is isomorphic to a twisted group 
ring $\mathbb E^t[\Ga]$ for the extension field $\mathbb E= 
\cent(\scL)^0$ of $\F$ and for the subgroup \/$\Gamma=$ 
supp$_{\cent(\scL)}$ of $\La$.
\end{lem}

\Proof   Any $\Delta$-graded Lie algebra  $\scL$ is perfect,  so 
$\cent(\scL)$ is commutative. As a $\mult(\scL)$-module, $\scL$ is 
generated by $\h$. Indeed, $[\scL_\al,\h]=\scL_\al$ for $0\ne \al 
\in \De$ and then $\scL_0 = \sum_{\al\in \De}\, [\scL_{-\al}, 
[\scL_\al, \h]]$ by (iv) in Definition \ref{rootgradE}. Since $\h$ 
is finite-dimensional, it follows from (\ref{fingen}) that $\cent 
(\scL) = \grcent(\scL)$ with respect to the $\La \oplus 
\scQ(\De)$-grading  (\ref{DeLamGrad}). However, by Proposition 
\ref{toral} (a), every $\chi$ has degree $0$ with respect to the 
$\scQ(\De)$-grading of $\scL$, which proves the first part of the 
lemma.  The last part now follows from Proposition \ref{gradprop}. 
\qed \psp 

\begin{defn}\label{deflietor}
A $(\De,\La)$-graded Lie algebra is called a {\it Lie torus} of type 
$(\De,\La)$, or simply a {\it Lie torus} if
\begin{itemize}  
\item[{\rm (i)}] If $\scL_\alpha^\lambda \neq 0$ for $\alpha,\lambda \neq 0$, then there exists
an $\mathfrak{sl}_2$-triple $(e,h,f) \in \scL_\alpha^\lambda  \times 
\scL_0^0 \times \scL_{-\alpha}^{-\lambda}$ such that 
\begin{itemize} 
\item[{(a)}]   $\scL_\al^\lambda = \F e$ and $\scL_{-\alpha}^{-\lambda} = \F f$, and
\item[{(b)}]   $[h, x] = \la \beta, \alpha^\vee \ra x$ for all $x \in \scL_{\beta}^\mu$,  \ 
$\beta \in \De \cup \{0\}$,  \ $\mu \in \Lambda$, and $h \in \h$;
\end{itemize}
\item[{\rm (ii)}]  The group $\La$ is generated by 
  supp$_{\scL} = \{\lambda \in \La \mid \scL^\lambda \neq 0\}$;
\item[{\rm (iii)}] If $\Delta = {\rm BC}_r$ for $r=1,2$, then 
$\Delta_\g \ne {\rm D}_r$. 
\end{itemize}  \end{defn} \psp

In (b), $\alpha^\vee$ is the usual {\it coroot}, and $ \la \beta, 
\alpha^\vee \ra$ comes from the usual inner product on the real span 
of the root system $\Delta$.   An $\mathfrak{sl}_2$-triple in (i)  
is assumed to satisfy the canonical relations, $[e,f] = h$, $[h,e] = 
2e$, and $[h,f] = -2f$. We point out that no condition is imposed  
on $\dim \scL^\lam_0$; however, when $\La$ is free of finite rank, 
one knows that $\dim \scL^\lam_0 \leq m$  for some positive integer 
$m$ that does not depend on $\lambda$  by  \cite[Thm.~5(a)]{N2}. In 
\cite{N2}, Neher considers Lie tori only for $\La = \ZZ^n$;   while 
in 
 \cite{Y4},  the above notion is referred to as a $\La$-torus, 
and the term Lie torus is reserved for the case $\La = \ZZ^n$.    
Here we consider  Lie tori for an arbitrary $\La$,  since the 
determination of the centroid  is exactly the  same as for $\La = 
\ZZ^n$.  \psp

A centreless Lie torus is graded-simple by \cite[Lem.~4.4]{Y4}.   It 
is assumed in \cite{Y4}  that $\Delta_\g$ is of type B$_r$ when 
$\Delta$ is of type BC$_r$;  however,  the same proof as given there 
works in our more general setting.   \psp
 
The following result provides a proof of \cite[Thm.~7(a)]{N2}: 
 \psp

\begin{prop} \label{centlt}
Let $\scL$ be a centreless Lie torus of type $(\De,\La)$. Then 
$\cent(\scL)$ is a twisted group ring $\F^t[\Gamma]$ for some 
subgroup $\Gamma$ of $\La$. In particular, 
\begin{enumerate}
 \item[\rm (i)] $\cent(\scL)^0 = \F\Id$, and 
 \item[\rm(ii)] if $\La \cong \ZZ^n$,  then $\cent(\scL)$ is isomorphic 
  to a Laurent polynomial ring.   
\end{enumerate}
\end{prop}

\Proof  We know from Lemma \ref{rootgradcen} and \cite[Lem.~3.4]{Y4}  
that $\cent(\scL) = \bigoplus_{\lam \in \La}\, \cent(\scL)^\lam$ is 
a commutative associative division-graded algebra.   Since 
$\cent(\scL)^\lam \to \scL^\lam_\al$ is injective for each $\lambda$ 
and $\alpha$ by Proposition \ref{gradprop}, it follows that all 
homogeneous spaces $\cent(\scL)^\lam$ are at most one-dimensional. 
In particular (i) holds.  It also follows that $\cent(\scL)$ is a 
twisted group ring for a subgroup $\Ga$ of $\La$.    When  $\La = 
\ZZ^n$, this subgroup is isomorphic to $\ZZ^r$ for some $0 \leq r 
\le n$. But any twisted group ring over $\ZZ^r$ is in fact a group 
ring,  and hence isomorphic to a Laurent polynomial ring, (compare  
\cite[Lem.~1.8]{BGKN}).   \qed \pse 

 \section{The centroid of an EALA and its core}\label{sec4}

We will prove in Corollary  \ref{ealacent} below that the core of a 
tame EALA  is central, and from this result,  the centroid of the 
EALA itself can easily be determined.  Since the core of an EALA is 
a central extension of a centreless Lie torus,  and since we know 
the centroid of centreless Lie tori by Proposition \ref{centlt}, it 
is natural to base our investigation on a general result which 
describes the behavior of the centroid  under a central extension, 
see Lemma \ref{xxx} below.   We start by recalling some facts about 
central extensions and $2$-cocycles, which also serves to establish 
our notation. \psp

Throughout we consider Lie algebras over an arbitrary field $\F$. We 
recall that a  {\it central extension} of a Lie algebra $\LL$ is a 
pair $(\KK,\pi)$ consisting of a Lie  algebra $\KK$ and a surjective 
Lie algebra homomorphism $\pi: \KK  \rightarrow \LL$ whose kernel 
lies in the center $Z(\KK)$ of $\KK$.    If $\KK$ is perfect, then 
$\KK$ is said to be a {\it cover} or {\it covering} of $\LL$.   In 
this case $\LL$ is necessarily perfect also.   A {\it homomorphism} 
(resp. {\it isomorphism\/}) from a central extension $f : \KK \to 
\LL$ to a central extension $f' : \KK' \to \LL$ is a Lie algebra 
homomorphism (resp. isomorphism) $g : \KK \to \KK'$ satisfying $f = 
f' \circ g$.  \psp

A central extension $\fru : \widehat \LL \to \LL$ is a {\it 
universal central extension} if there exists a unique homomorphism 
from $\widehat \LL$ to any other central extension $\KK$ of $\LL$.  
This universal property implies that any two universal central 
extensions of $\LL$ are isomorphic as central extensions. A Lie 
algebra $\LL$ has a universal central extension if and only if $\LL$ 
is perfect. In this case, the universal central extension $\widehat 
\LL$ is perfect,  and $\widehat \LL$ is a covering of every covering 
of $\LL$. \psp

Examples of central extensions can be constructed in terms of {\it 
$2$-cocycles} which are bilinear maps $\sig: \LL \times \LL \to C$ 
into some vector space $C$ satisfying
\begin{equation}\label{(1)}
\si(x,x) = 0 \quad\hbox{and}\quad \si([x,y],z) + \si([y,z],x) + 
\si([z,x],y) = 0
\end{equation}
for all $x,y,z\in \LL$.  Given such a $2$-cocycle, the vector space 
$E = \LL \oplus C$ becomes a Lie algebra with product
\begin{equation} \label{centbe}
[x \oplus c \,,\, y \oplus c']_E =  [x,y]  \oplus \si(x,y).
\end{equation}
 (Here we are using the notation $x \oplus c$ to designate
that $x \in \LL$ and $c \in C$.)  This is  a central extension of 
$\LL$, which we denote $E(\LL,\sig)$, with respect to the projection 
map $E(\LL,\sig) \to \LL, \quad  x\oplus  c \mapsto x$. Since we are 
considering Lie algebras over a field, every central extension $f: 
\KK \to \LL$ is isomorphic as central extension to some 
$E(\LL,\sig)$, see e.g. \cite[1.9]{MP}. \psp

\begin{lem}\label{xxx}  Let $\pi : \KK  \to
\scL$ be a central extension of the Lie  algebra $\LL$ written in 
the form {\rm (\ref{centbe})}, and suppose that $Z(\LL) = 0$. Then $ 
Z({\KK}) = C$.  Moreover, $\Psi\in \End_\F (\KK)$ lies in the 
centroid $\cent(\KK)$ if and only if there exist $\chi \in 
\cent(\LL)$, $\psi \in {\rm Hom}_\F(\scL,C)$ and $\eta \in 
\End_\F(C)$ such that
\begin{eqnarray}
\Psi(x \oplus  c) &=& \chi(x) \oplus  \bigl(\psi(x) + \eta(c)\bigr) 
\quad\hbox{and} \label{xxx2} \\ \sigma(x,\chi(y)) &=& \psi([x,y]) + 
\eta\bigl(\sigma(x,y)\bigr) \label{xxx3}
\end{eqnarray}
 for all $x,y\in \scL$ and $c\in C$.    In this case, $\sigma(x, \chi(y)) =  \si(\chi(x),y)$.
 \end{lem}
\medbreak

\Proof  If $x \oplus c \in Z(\KK)$,  then (\ref{centbe}) implies 
that $x\in Z(\LL)=0$, hence $Z(\KK) \subseteq C$. The other 
inclusion is obvious.

Now assume $\Psi \in \cent(\KK)$. By Lemma~\ref{oblem}\,(a), $\Psi$ 
leaves $Z(\KK)=C$ invariant, hence has the form (\ref{xxx2}) for 
$\psi$, $\eta$ as in the statement of the lemma and some $\chi \in 
\End_\F(\scL)$. Since $\pi_{\cent} (\Psi) = \chi$, it follows from 
Lemma~\ref{oblem}\,(a) that $\chi\in \cent(\scl)$ (this is also 
immediate from the computation below).   It now remains to 
characterize when a map of the form (\ref{xxx2}) belongs to the 
centroid of $\KK$. For $x,y \in \scl$ and $c,c'\in C$ we have
 \begin{eqnarray*}
&& \Psi\left ([x \oplus  c, y \oplus c']_{\KK}\right)  =  \Psi\bigl([x,y] \oplus \sigma(x,y)\bigr) \\
&&  \hskip 1.3 truein =
\chi([x,y]) \oplus \Big(  \psi([x,y]) + \eta\bigl (\si(x,y) \bigr)\Big), \\
& & [x\oplus c, \Psi(y \oplus c')]_\KK  =  [x,\chi(y)] \oplus 
\si(x,\chi(y)).   \end{eqnarray*}  Combined they show that $\Psi \in 
\cent(\KK)$ if and only if  $\chi \in \cent(\LL)$ and (\ref{xxx3}) 
holds. \qed \psp

Next we will construct a special class of $2$-cocycles for Lie 
algebras with a nondegenerate invariant bilinear form and describe 
the centroid of the corresponding central extension in  Proposition 
\ref{centprop}.   Later this will be applied to determine the 
centroid of the core of an EALA. \psp

Let $\Lambda$ be an abelian group.  We say $\pi : \KK  \to \LL$ is a 
{\it $\Lambda$-graded central extension} if both $\LL$ and $\KK$ are 
graded by $\Lambda$, say $\LL = \bigoplus_{\lambda \in \Lambda} \, 
\LL^\lambda$ and $\KK = \bigoplus_{\lambda \in \Lambda}\, 
\KK^\lambda$, and if $\pi$ is homogeneous of degree $0$, i.e., 
$\pi(\KK^\lambda) \subseteq \LL^\lambda$.    Every $\Lambda$-graded 
central extension $\pi: \KK \to \LL$ is isomorphic to a central 
extension $E(\LL,\sig)$ where $\sig : \LL \times \LL \to C$ is a 
{\it $\Lambda$-graded $2$-cocycle}, i.e., $C= \bigoplus_{\lambda \in 
\Lambda} C^\lambda$ is a $\Lambda$-graded vector space and $\sig$ is 
a $2$-cocycle satisfying $\sig(\LL^\lambda, \LL^\mu) \subseteq 
C^{\lambda + \mu}$ for all $\lambda,\mu \in \Lambda$.   When $\KK$ 
is perfect,  such a $\Lambda$-graded central extension  $\pi : \KK  
\to \LL$  is called a  {\it $\Lambda$-covering}.    In particular, 
if $\LL$ is perfect, the universal central extension $\widehat \LL$ 
is a $\Lambda$-covering of $\LL$. \medbreak

Let $\LL= \bigoplus_{\lambda  \in \Lambda}\, \LL^\lambda$  be a 
$\Lambda$-graded Lie algebra over $\F$.   We denote by 
$\Hom_\ZZ(\Lambda,\F)$ the $\F$-vector space of additive maps 
$\theta : \Lambda  \to \F$.    For any $\theta \in 
\Hom_\ZZ(\Lambda,\F)$,   the corresponding  {\it degree derivation} 
$\pa_\theta$ of $\LL$  is defined  by
\begin{equation}\label{degder} \pa_\Th (x) = \Th(\lambda) x \qquad (x  \in \LL^\lambda), 
\end{equation} 
There is a linear map
\begin{equation*}
  \Hom_\ZZ(\Lambda, \F) \to \scd : = \{ \pa_\Th \mid  \Th \in \Hom_\ZZ(\Lambda, \F)\} \qquad  \Th \mapsto \pa_\Th,   \end{equation*}
into the space $\scd$ of degree derivations,  which is an 
isomorphism if $\Lambda$ is spanned by the support of $\LL$, i.e., 
if
 \begin{equation} \Lambda = \hbox{\rm span}_\ZZ \{ \lambda \in \Lambda \mid \LL^\lambda
\ne 0 \}.  \label{supp}
 \end{equation}
Indeed, if $\pa_\Th = 0$ then $\Th(\lambda) = 0$ for all $\lambda 
\in \Lambda$ with $\LL^\lambda \ne 0$ whence $\Th=0$. We note that 
(\ref{supp}) is essentially a notational convenience;   if it is not 
fulfilled,  one can always replace $\Lambda$ by the subgroup 
generated by the support of $\LL$. Assuming (\ref{supp}), we have a 
well-defined linear map $$\ev : \Lambda \to {\mathcal D}^*, \quad  
\lambda \mapsto \ev(\lambda)$$ into the dual space $\mathcal D^*$ of 
$\mathcal D$ given by
\begin{equation}
\ev(\lambda)(\pa_\Th) = \Th(\lambda), \quad
   \lambda \in \Lambda,\, \, \pa_\Th \in {\mathcal D}. \label{ev}\end{equation}
We denote by ${\rm Der}( \LL)^\lambda$,  the vector space of 
$\F$-linear derivations of $\LL$ of degree $\lambda$,  and we set
$$  {\rm gr}\Der( \LL) = \bigoplus_{\lambda \in \Lambda}\, \Der( \LL)^\lambda,  
$$
which is obviously a graded subalgebra of  $\Der(\LL)$. It is 
well-known (see for example \cite[Prop.~1]{F}) that ${\rm grDer}( 
\LL) = \Der( \LL)$ if $\LL$ is finitely generated as Lie algebra. 
\psp

Let $(\,.\mid.\,)$ be an invariant $\Lambda$-graded bilinear form on 
$\LL$, i.e., $(\LL^\lambda \,|\, \LL^\mu) = 0$ for all $\lambda +\mu  
\ne 0$. We denote by SDer$(\LL) =$ SDer$_\F(\LL)$ the subalgebra of 
$\Der( \LL)$ consisting of skew derivations $\partial$ of $\LL$:  
$(\partial x \,|\, y) = - (x \,|\, \partial y)$ for all $x,y \in 
\LL$.     It is easily seen that $\pa \in {\rm grDer}(\LL)$ is skew 
if and only if every homogeneous component of $\pa$ is,  so that
$$ \SDer( \LL) := \rm{SDer}(\LL) \cap {\rm grDer}(\LL) =  \bigoplus_{\lambda \in \Lambda} \,  \hbox{\rm SDer}( \LL)^\lambda.
$$
where ${\rm SDer}( \LL)^\lambda = {\rm SDer}(\LL) \cap 
\Der(\LL)^\lambda$.   Moreover
 $$   {\mathcal D} \subseteq  {\rm SDer}(\LL)^0 \quad\hbox{and}\quad
    {\rm IDer}(\LL) \subseteq \SDer(\LL)
$$
where ${\rm IDer} (\LL) = {\rm IDer}_\F(\LL) = \{ \ad x \mid  x\in 
\LL \}$ denotes the ideal of inner derivations of $\LL$. \psp

Let ${\mathcal S}=\bigoplus_{\lambda \in \Lambda} {\mathcal 
S}^\lambda$ be a graded subspace of $\SDer(\LL)$, and let ${\mathcal 
S}^{{\rm gr}*}$ be the graded dual space. Thus,
\begin{equation}\label{eq:grdual} {\mathcal S}^{{\rm gr}*} = \bigoplus_{\lambda \in \Lambda} ({\mathcal S}^{{\rm gr}*})^\lambda \, ,\quad \hbox{\rm where} \quad 
        ({\mathcal S}^{{\rm gr}*})^\lambda= ({\mathcal S}^{-\lambda})^*.
\end{equation}
We may assume  $({\mathcal S}^{-\lambda})^* \subseteq {\mathcal 
S}^*$ by defining $f\mid_{{\mathcal S}^\mu} = 0$ for $f\in 
({\mathcal S}^{-\lambda})^*$ and $\mu \ne -\lambda$.    Then it is 
easy to verify that
\begin{equation}
\label{Scocy} \si_{\mathcal S} : \LL \times \LL \to {\mathcal 
S}^{{\rm gr}*}, \quad \si_{\mathcal S}(x,y)(d)= (d(x)\,|\, y)
\end{equation}
for $x,y\in \LL$ and $d\in {\mathcal S}$ is a $\Lambda$-graded 
$2$-cocycle. Thus  $E(\LL,\sig_{\mathcal S}) = \LL \oplus {\mathcal 
S}^{{\rm gr}*}$ with product $[x \oplus c, y\oplus c']_E =[x,y] 
\oplus \sigma_{{\mathcal S}}(x,y)$ for all $x,y \in \LL, c,c' \in  
{\mathcal S}^{{\rm gr}*}$ is a $\Lambda$-graded central extension of 
$\LL$. \psp

\begin{prop}\label{centprop}
Let $\LL=\bigoplus_{\lambda \in \Lambda}\, \LL^\lambda$ be a perfect 
$\Lambda$-graded Lie algebra with a nondegenerate invariant graded 
bilinear form such that  {\rm (\ref{supp})} holds. Let ${\mathcal 
S}\subseteq \SDer( \LL)$ be a $\Lambda$-graded subspace such that
\begin{equation}\label{evd}
    \ev_{\mathcal S} : \Lambda  \to ({\mathcal D}\cap {\mathcal S})^*, \quad \lambda  \mapsto
    \ev(\lambda)\,\big|_{{\mathcal D} \cap {\mathcal S}}
\end{equation}
is injective, where the evaluation map $\ev$ is as in {\rm 
(\ref{ev})}. Let $\KK=E(\LL,\sig_{\mathcal S}) = \LL \oplus C$ where 
$C = {\mathcal S}^{{\rm gr}*}$.   \begin{itemize}
\item[{\rm (i)}] Suppose $\Psi \in \cent(\KK)$ is homogeneous of degree 
$\lambda$. Then there exists $\chi \in \cent(\LL)$, $\psi \in 
\hbox{\rm Hom}_\F(\LL,C)$,  and $\eta \in \hbox{\rm End}_\F(C)$ all 
of degree $\lambda$ such that
\begin{itemize}
\item[{\rm (a)}]  $\Psi\bigl(x \oplus c) =\chi(x)  \oplus
\bigl(\psi(x) \oplus \eta(c) \bigr)$   for all  $x \in \LL, c \in 
C$. \item[{\rm (b)}]   $\chi = 0$ if  $\lambda \neq  0$.
\end{itemize}
\item[{\rm (ii)}]   If $\KK$ is perfect,  then $\cent(\KK)^\lambda = 
0$ for  all $\lambda  \neq 0$. In particular, if $\scK$ is a Lie 
torus, then $\KK$ is central.   \end{itemize}
\end{prop}
\medbreak

\Proof  (i) Observe first that $Z(\LL) = 0$.   This follows from the 
computation $([x,y] \mid z) = (x \mid [y,z]) = 0$ for all $x,y \in 
\LL$, $z \in Z(\LL)$ and the fact that $\LL$ is perfect and the form 
is nondegenerate.

Now suppose $\Psi \in \cent(\KK)$ has degree $\lambda$,  and apply 
Lemma \ref{xxx} to conclude that $\Psi(x \oplus c) = \chi(x) \oplus 
\bigl(\psi(x) + \eta(c)\bigr)$ where $\chi \in \cent(\LL)$, $\psi 
\in \hbox{\rm Hom}_\F(\LL,C)$, and $\eta \in \hbox{\rm End}_\F(C)$,  
and all have degree $\lambda$.   For $x^\mu \in \LL^\mu$, $y \in 
\LL$,  and $\pa_\Th \in {\mathcal D} \cap {\mathcal S}$ we have by 
(\ref{xxx3}),
\begin{eqnarray*}
\sig_{\mathcal S}(\chi(x^u)\,,\,y)(\pa_\Th) &=& (\pa_\Th \chi(x^\mu) 
\,|\, y )
  = \Th(\mu +\lambda) ( \chi(x^\mu) \,|\, y) \\
  &=& \sig_{\mathcal S}(x^\mu, \chi(y))(\pa_\Th) = (\pa_\Th(x^\mu) \,|\,\chi (y))
  = \Th(\mu)(x^\mu \,|\, \chi(y))\\
  & =& \Th(\mu)(\chi(x^\mu) \,|\,y),
\end{eqnarray*}
where in the last equality we used Lemma~\ref{easy}\,(f). Hence 
$\Th(\lambda)(\chi(x^\mu) \,|\, y) = 0$ for all $x^\mu \in \LL^\mu$ 
and $\mu \in \Lambda$.    Suppose $\lambda \neq 0$.   Then there 
exists a $\theta  \in {\mathcal D} \cap {\mathcal S}$  such that 
$\theta(\lambda) \neq 0$.  From this we see $(\chi(x^\mu), y) = 0$ 
and nondegeneracy forces $\chi = 0$. Thus, (b) holds. 

(ii) In the preceding paragraph we have shown that  $Z(\LL) = 0$, 
and $\chi = 0$ whenever $\lambda \neq 0$.   Thus if $\pi: \KK 
\rightarrow \LL$ is the cover map and  $\KK$ is perfect,  it follows 
from Lemma \ref{oblem}\,(c) that $\Psi = 0$. If $\KK$ is a Lie 
torus, we can apply Proposition~\ref{centlt} to conclude that $\KK$ 
is central. \qed

\psp As an application of Proposition \ref{centprop},  we can now 
determine the centroid of a tame EALA $\scE$  and of  its core 
$\scK$, which is the ideal of $\mathcal E$ generated by the root 
spaces corresponding to the nonisotropic roots.  The reader is 
referred to \cite[Ch.~I]{AABGP} for the precise definition of a tame 
EALA  over $\F=\CC$ and to \cite{N3} for   arbitrary fields $\F$ of 
characteristic $0$.      An EALA has an invariant nondegenerate 
bilinear form,  which when restricted to the core $\KK$  is 
$\Lambda$-graded for $\Lambda= \ZZ^n$ for some $n \geq 0$. Tameness 
says that the ideal $\scK$ satisfies  $\scC_{\scE}(\scK) = Z(\scK)$.  

\begin{cor}\label{ealacent} Let  $\scE$ be a tame extended affine 
Lie algebra, let $\KK$ be its core and set $\scD=\scE/\scK$. Then 
$\KK$ is central and 
$$\cent(\scE) = \F\Id \oplus \scV(\KK) 
  \cong \F\Id \oplus \Hom_\scD (\scD, Z(\KK)).
  $$
\end{cor}

\Proof   It is known that $\KK$ is a Lie torus with $\La \cong 
\ZZ^n$ (see \cite[Cor.~7.3]{Y5}  for  $\F=\CC$ or 
\cite[Prop.~3(a)]{N3} for arbitrary $\F$). Moreover, by 
\cite[Thm.~6]{N3}, $\KK$ is  obtained from the centreless Lie torus 
$\scL = \scK /Z(\scK)$ by the construction of Proposition 
\ref{centprop}. Thus $\scK$ is central by part (ii) of that 
proposition.  Since $\scK$ is a $\cent(\scE)$-invariant ideal of 
$\scE$, \ $\cent(\scE) = \F\Id \oplus \scV(\scK)$ follows from 
(\ref{smacent}). Finally, $\scV(\scK) \cong \Hom_\scD (\scD, 
Z(\KK))$ by Lemma \ref{easy}\,(b) and the fact that 
$\Ann_{\scE}(\scK) = \scC_\scE (\scK) = Z(\scK)$ because of 
tameness. \qed

\medbreak

\begin{exam}\label{examEALA}
{\rm Finite-dimensional split simple Lie algebras are examples of 
tame EALAs.   In this case $\scE=\scK$ and $Z(\scK)=0$, so the 
result above simply says that $\scE$ is central -- which is of 
course well-known.  

Another class of examples of tame EALAs are the affine Lie algebras 
(see \cite{ABGP}).   In this case, $\scK=\scE^{(1)}$ and $\scD$ and 
$Z(\scK)$ are both one-dimensional, so that our result recovers 
Corollary~\ref{centKM}. 
 
However there are many other examples of extended affine Lie 
algebras besides the two just mentioned, e.g. toroidal algebras 
extended by some derivations. }\end{exam} \pse

\section{Centroids of Lie Algebras Graded by Finite Root Systems}
\label{rootgrad}

In this section, we will describe the centroid of Lie algebras 
graded by finite root systems. The case of reduced root systems will 
be treated in Subsection \ref{red}, while the nonreduced case will 
be done in \ref{nonred}.   As a prelude to that, we begin with a 
general result about Lie algebras $\LL$ which are completely 
reducible relative to the adjoint action of a subalgebra $\g$. By 
gathering together isomorphic summands, we may assume that such Lie 
algebras are written in the form
$$\LL = \bigoplus_{k}  \left(V_k \ot A_k\right),$$
where the $V_k$ are nonisomorphic irreducible $\g$-modules;  the 
subspace $A_k$ indexes  the copies of $V_k$; \ and the $\g$-action 
is given by $x.(v_k \ot a_k) = [x,v_k \ot a_k)]  = x.v_k \ot a_k$ 
for $x \in \g$, $v_k \in V_k$, $a_k \in A_k$.  \psp

\begin{lem}\label{lemcr}  Assume  $\LL$ is a Lie algebra which is completely reducible
relative to the adjoint action of a subalgebra $\g$, and let $\LL = 
\bigoplus_{k}  \left(V_k \ot A_k\right)$ be its decomposition 
relative to $\g$. Assume $\End_{\g}(V_k) = \F \Id$ for each 
irreducible $\g$-module $V_k$. If $\Psi \in \cent(\LL)$, then there 
exist transformations $\psi_k: A_k \to A_k$ such that $\Psi(v_k \ot 
a_k) = v_k \ot \psi_k(a_k)$ for all $v_k \in V_k$, $a_k \in A_k$.   
\end{lem}

\proof   Assume $\{a_k^i\mid i \in \mathfrak{I}_k\}$ is a basis for 
$A_k$, and let $\pi_k^i$ denote the projection of $\LL$ onto the 
summand  $V_k \ot a_k^i$.   Fix $j \in \mathfrak {I}_k$.  Then for 
any $\Psi \in {\rm Cent}(\LL)$ and any $i \in \mathfrak {I}_k$, we 
have $(\pi_k^i \circ \Psi):  V_k  \ot a_k^j \to  V_k \ot a_k^i$ is a 
$\g$-module homomorphism.  Thus, it determines an element of 
$\End_{\g}(V_k) = \F \Id$, and there exists a scalar $\xi_{i,j} \in 
\F$, so that $(\pi_k^i \circ \Psi)(v_k \ot a_k^j) = \xi_{i,j} v_k 
\ot a_k^i$ for all $v_k \in V_k$.  When $\ell \neq k$,   $(\pi_k^i 
\circ \Psi):  V_\ell  \ot a_\ell^j \to  V_k \ot a_k^i$ is a 
$\g$-module homomorphism determining an element of 
$\Hom_{\g}(V_\ell,V_k)$.  Such a homomorphism must be the zero map 
since $V_k$ and $V_\ell$ are irreducible and nonisomorphic.    
Consequently, $\Psi(v_k \ot a_k^j) \in V_k \ot A_k$ for all $v_k \in 
V_k$, and
$$\Psi(v_k \ot a_k^j) = \sum_{i \in \mathfrak {I}_k} \xi_{i,j} v_k \ot a_k^i
= v_k \ot \Bigg( \sum_{i \in \mathfrak {I}_k} \xi_{i,j} a_k^i 
\Bigg).$$

\noindent  Define $\psi_k(a_k^j) = \sum_{i \in \mathfrak {I}_k} 
\xi_{i,j} a_k^i$ for each $j \in \mathfrak{I}_k$ and extend this 
linearly to all of $A_k$. Then
$$\Psi(v_k \ot a_k) = v_k \ot \psi_k(a_k) \qquad \hbox{\rm for
all} \ \ v_k \in V_k, \ a_k \in A_k.   \qed$$  \pse

 \subsection{Lie algebras graded by reduced root systems}\label{red}
 
Suppose that $\LL$ is a Lie algebra graded by the reduced root 
system $\Delta$ as in Definition \ref{rootgradE}.   Then $\LL$ is 
completely reducible relative to the adjoint action of the grading 
subalgebra $\g$, and by results in (\cite{BeM}, \cite{BZ}, 
\cite{N1}) we know that
$$\LL \cong (\g \otimes A) \oplus (W \otimes B) \oplus D,$$
\noindent  where the following hold:

 \begin{itemize}
\item[{\rm (1)}]  $W$ is the irreducible $\g$-module with highest
weight the highest short root of $\g$;   thus $W$ and $B$ are zero 
when $\Delta$ is simply laced. \item[{\rm (2)}]   The sum $\mathfrak 
a = A \oplus B$ is a unital algebra called the {\em coordinate 
algebra} of $\LL$.   In all cases except for type C$_2$, \ 
$\mathfrak a$ is an associative, alternative, or Jordan algebra 
depending on $\Delta$.   The unit element $1$ of $\mathfrak a$ lives 
in $A$, and $\g$ is identified with $\g \ot 1$.  \item[{\rm (3)}]  
$D$ is the sum of the trivial one-dimensional $\g$-modules. 
Moreover, $D$ can be identified with a quotient space, $D = \langle 
\mathfrak a, \mathfrak a\rangle = \langle A,A\rangle + \langle 
B,B\rangle$, of the skew-dihedral homology of $\mathfrak a$, and $D$ 
acts by inner derivations on $\mathfrak a$.   Thus, $\langle \al, 
\be \rangle (\gamma) = D_{\al,\be}(\gamma)$ for all $\al,\be,\gamma 
\in \mathfrak a$, where $D_{\al,\be}$ is the inner derivation 
determined by $\al,\be$.  \item[{\rm (4)}] the multiplication in 
$\LL$ is given in terms of the product on $\mathfrak a$ as follows 
(note here we do not use $\oplus$ to separate summands to simplify 
the expressions):
 \end{itemize}
 \begin{equation}\label{firsttype} \big(\Delta = \text{\rm A}_1, \ \text{\rm
B}_r, (r \geq 3), \ \text{\rm D}_r, (r \geq 4), \ \text{\rm E}_6,\ 
\text{\rm E}_7, \ \text{\rm E}_8,\ \text{\rm F}_4, \ \text{\rm G}_2 
\big)\end{equation} \vskip -.3truein
\begin{eqnarray*} &&[x \ot a, y \ot a^\prime] = [x, y]\ot aa^\prime + (x|y)\la
a,a^\prime\ra\\
&&[d, x \ot a] = x \ot da = - [x \ot a, d]  \\
&&[x \ot a, u \ot b] =  xu \ot a b  = -[u \ot b, x \ot a]  \\
&&[d, u \ot b] = u \ot db = - [u \ot b, d] \\
&&[u \ot b, v \ot b^\prime] = \partial_{u,v} \ot (b,b^\prime) + 
(u\ast v)\ot(b\ast b^\prime) +(u |v)\la b,b^\prime\ra,
\end{eqnarray*}

\begin{equation}\label{secondtype} \big (\Delta =\text{\rm A}_r, (r \geq 2), \;
\text{\rm C}_r, (r \geq 2) \big )\end{equation}

\vskip -.3 truein

\begin{eqnarray*} && [x\ot a, y \ot a^\prime] = [x,y]\ot \frac
{1}{2}(aa^\prime + a^\prime a) + (x \circ y)\ot \frac {1}{2}(a
a^\prime -a^\prime a) + (x|y)\la a,a^\prime\ra, \nonumber \\
&& [d, x \ot  a ] = x\otimes da =  -[x \ot
a, d]  \\
&& [x \ot a, u \ot b] =  (x \circ u )\ot \frac {1}{2}(a b -b a) + 
[x,u] \ot \frac {1}{2}(ab +
ba) = -[u \ot b, x \ot a]   \\
&& [u \ot b, v \ot b^\prime] =  [u,v] \ot \frac {1}{2}(bb^\prime + 
b^\prime b) + (u \circ v) \ot \frac {1}{2}(bb^\prime -
b^\prime b) + (u|v)\la b,b^\prime\ra  \\
&& [d,u \ot b] = u \ot db = -[u \ot b,d]
\end{eqnarray*}
\begin{itemize}
\item[{}]   for all $a,a' \in A, \, b,b' \in B$, $x,y \in \g$,
$u,v \in W$, and $d \in D$.  \item[{\rm (5)}]  In (4),  $(\,|\,)$ 
denotes the Killing form when applied  to $\g$, and the unique 
$\g$-invariant bilinear form when applied to $W$.  The maps 
$\partial  \in \Hom_\g(W\ot W,\g)$ and $\ast \in \Hom_{\g}(W \ot W, 
W)$ in (\ref{firsttype}) are unique up to scalars.    In 
(\ref{firsttype}),  the product in the coordinate algebra $\mathfrak 
a = A \oplus B$  is given by $(a + b)(a' + b') = \Big(aa' + 
(b,b')\Big) +  \Big(ab' + a'b + b \ast b' \Big)$. The algebra $\g$ 
in (\ref{secondtype}) can be realized as a matrix Lie algebra 
$\mathfrak{sl}_{r+1}$ or $\mathfrak{sp}_{2r}$ respectively. So for 
any two matrices $w,z$, we have  $[w,z] = wz-zw$ and $w \circ z = 
wz+zw - (2/n) \mathfrak{tr}(wz)$,  where $n = r+1$ or $2r$ and 
$\mathfrak {tr}$ denotes the trace.
\end{itemize}
\psp

\begin{prop}\label{prop1rg}  Assume $\LL \cong (\g \otimes A) \oplus (W \otimes B) \oplus D$
 is a Lie algebra graded by a finite reduced root system. Then $\Psi 
 \in \cent(\LL)$ if and only if there exist maps $\psi_{\mathfrak a} \in \cent(\mathfrak a)$ and $\psi_{_D} \in
\cent(D)$  such that   $\psi_{\mathfrak a}(A) \subseteq A$, 
$\psi_{\mathfrak a}(B) \subseteq B$, and
\begin{eqnarray}\Psi\Big ( (x \ot a) + (u \ot b)+ d \Big) 
  &=& \bigl(x \ot \psi_{\mathfrak a}(a)\bigr)
+ \bigl(u \ot \psi_{\mathfrak a}(b)\bigr) + \psi_{_D}(d) \qquad  \\
\psi_{_D}\big (\la \alpha, \alpha'\ra \big ) &=& \la \alpha, 
\psi_{\mathfrak a}(\alpha') \ra = \la
\psi_{\mathfrak a}(\alpha), \alpha' \ra \label{vard} \\
\big(\psi_{\mathfrak a}\circ d\big)(\alpha) &=& \big(d \circ
\psi_{\mathfrak a}\big)(\alpha) \\
\psi_{_D}(d)(\alpha) &=& \big(\psi_{\mathfrak a} \circ d \big) 
(\alpha)
\end{eqnarray}
\noindent  for all $\alpha, \alpha' \in \mathfrak a$ and $d \in D$.    
\end{prop}

\proof  \ Applying Lemma  \ref{lemcr},  we see that corresponding to 
$\Psi \in \cent(\LL)$ are maps $\psi_{_A} \in \End_\F(A)$, 
$\psi_{_B} \in \End_\F(B)$, and $\psi_{_D} \in \End_\F(D)$  such 
that $\Psi(x \ot a) = x \ot \psi_{_A}(a)$,  $\Psi(u \ot b) = u \ot 
\psi_{_B}(b)$ and $\Psi(d) = \psi_{_D}(d)$ for $x \ot a \in \g \ot 
A$, $u \ot b \in W \ot B$, and $d \in D$.    Set $\psi_{\mathfrak 
a}(a \oplus b) = \psi_{_A}(a) + \psi_{_B}(b)$ and observe that 
$\psi_{\mathfrak a}(A) \subseteq A$ and $\psi_{\mathfrak a}(B) 
\subseteq B$ clearly hold. 

Now suppose $ x \ot a, y \ot a' \in \g \ot A$, and consider $\Psi([x 
\ot a, y \ot a'])$.  When $\Delta$ is as in (\ref{firsttype}), then
\begin{eqnarray*}
\Psi([x \ot a, y \ot a']) &=& [x\ot a,\Psi(y \ot a')] \qquad  \Longleftrightarrow\\
{[x,y] \ot \psi_{\mathfrak a}(aa') + (x|y)\psi_{_D}\big(\la 
a,a'\ra}\big) &=& [x \ot
a, y \ot \psi_{\mathfrak a}(a')] \\
&=& [x,y] \ot a \psi_{\mathfrak a}(a') + (x|y) \la a, 
\psi_{\mathfrak a}(a')\ra. \end{eqnarray*}

\noindent  Equating components shows that $\psi_{\mathfrak a}(aa') = 
a \psi_{\mathfrak a}(a')$,  and  $\psi_{_D}\big (\la a, a'\ra \big ) 
= \la a, \psi_{\mathfrak a}(a') \ra$ hold for all $a,a' \in A$.   
Similarly, using $\Psi([x \ot a, y \ot a']) = [\Psi(x \ot a), y \ot 
a']$, we obtain $\psi_{\mathfrak a}(aa') = \psi_{\mathfrak a}(a)a'$ 
and the second equality in (\ref{vard}) for all $a,a' \in A$. 

Now when $\Delta$ is as in (\ref{secondtype}), then
\begin{eqnarray}\label{AC} \Psi([x \ot a, y \ot a']) &=& [x,y]\ot \frac
{1}{2}\psi_{\mathfrak a}(aa^\prime + a^\prime a) \\
&& \ \  + (x \circ y)\ot \frac {1}{2}\psi_{\mathfrak a}(a a^\prime
-a^\prime a) + (x|y) \psi_{_D}\big(\la a,a^\prime\ra\big) \ \  \hbox{\rm while} \nonumber \\
{[x \ot a, \Psi(y \ot a')]} &=& {[x \ot a, y \ot \psi_{\mathfrak 
a}(a')]} \nonumber \\ &=& [x,y]\ot \frac {1}{2}\big(a\psi_{\mathfrak 
a}(a^\prime) + \psi_{\mathfrak a}
(a^\prime) a) \nonumber \\
&& \ \ + (x \circ y)\ot \frac {1}{2}\big(a \psi_{\mathfrak 
a}(a^\prime) - \psi_{\mathfrak a}(a^\prime) a\big)  + (x|y) \la 
a,\psi_{\mathfrak a}(a^\prime)\ra. \nonumber
\end{eqnarray}

\noindent  In particular, if $\Delta = \text{\rm A}_r$ for $r \geq 
2$, we may set $x = e_{1,1}-e_{2,2} = y$ (matrix units) and get 
$[x,y] = 0$, but $x \circ y \neq 0$.  Then equating components in 
these expressions, we obtain that
$$\psi_{\mathfrak a}Y(aa^\prime - a^\prime a) =
a \psi_{\mathfrak a}(a^\prime) - \psi_{\mathfrak a}(a^\prime) a$$

\noindent  holds for all $a,a' \in A$, (as does the first equality 
in (\ref{vard})).  Then putting that relation back in, we determine 
that
$$\psi_{\mathfrak a}(aa^\prime + a^\prime a) =
a \psi_{\mathfrak a}(a^\prime) + \psi_{\mathfrak a}(a^\prime) a$$

\noindent  for all $a,a' \in A$.  Combining these gives 
$\psi_{\mathfrak a}(aa') = a\psi_{\mathfrak a}(a')$ for all $a,a' 
\in A$.

When $\Delta = \text{\rm C}_r$ for $r \geq 2$, then the summands on 
the right side of (\ref{AC}) lie in different components, so 
equating them gives the same information as obtained in the A case.

Applying $\Psi$ to all the various other products in 
(\ref{firsttype}) and (\ref{secondtype}) and arguing similarly will 
complete the proof of the proposition.    \qed

\psp
 
For the centroidal transformation  $\psi_{\mathfrak a} \in 
\cent(\mathfrak a)$ coming from an element $\Psi \in \cent(\LL)$ as 
in Proposition \ref{prop1rg}, it  follows that $\psi_{\mathfrak 
a}(1)  \in A \cap Z(\mathfrak a)$, as $\psi_{\mathfrak a}$ preserves 
the space $A$ and the unit element of the coordinate algebra 
$\mathfrak a$ belongs to $A$.  \psp

The inner derivations of the coordinate algebra $\mathfrak a$ 
involve certain expressions in the left multiplication and right 
multiplication operators  which can be found in \cite[(2.41)]{ABG1}.
 For any $\Psi \in \cent(\LL)$, the associated  transformation
$\psi_{\mathfrak a}$ belongs to $\cent(\mathfrak a)$, so it commutes 
with the left and right multiplication operators of $\mathfrak a$.  
It also commutes with the involution $\sigma$ on $\mathfrak a$ when 
$\Delta$ is of type C$_r$, as $\psi_{\mathfrak a}$ preserves the 
spaces $A$ and $B$, which are the symmetric elements and 
skew-symmetric elements respectively relative to $\sigma$.     Thus, 
$\psi_{\mathfrak a}$ commutes with the inner derivation  
$D_{\al,\be}$ for all $\al,\be \in \mathfrak a$, and
\begin{equation}\label {eq:psicommute} \psi_{\mathfrak a} \circ D_{\al,\be} = D_{\psi_{\mathfrak a}(\al),\be} = D_{\al,\psi_{\mathfrak a}(\be)}. \end{equation}
As
\begin{equation}\la \al, \be \ra (\gamma) = D_{\al,\be}(\gamma), \end{equation}
\noindent  we have

\begin{equation}\Big(\psi_{\mathfrak a} \circ \la \al, \be \ra\Big)(\gamma)
= \la \psi_{\mathfrak a}(\al),\be\ra(\gamma) = \la 
\al,\psi_{\mathfrak a}(\be)\ra (\gamma) = \psi_{_D}(\la \al, 
\be\ra)(\gamma), \end{equation}   \noindent  so that for $\mathfrak 
z : = \psi_{\mathfrak a}(1)$, 

$${\mathfrak z} \la \al, \be \ra(\gamma) = \la \mathfrak z \al, \be \ra (\gamma) =
\la \al, {\mathfrak z} \be \ra (\gamma) = \psi_{_D}(\la \al, \be 
\ra)(\gamma)$$ \noindent  for all $\al,\be,\gamma \in \mathfrak a$.   
Combining this with Proposition \ref{prop1rg}, we obtain the 
following:

\psp

\begin{cor} \label{oneway} Assume $\LL$ is a $\Delta$-graded Lie algebra as in
Proposition \ref{prop1rg}, and let $\Psi \in\cent(\LL)$.   Then 
there exists an element ${\mathfrak z} \ (= \psi_{\mathfrak a}(1)\in 
A$) in the center $Z(\mathfrak a)$ of the coordinate algebra 
$\mathfrak a$ of $\LL$ such that
\begin{eqnarray}\label{cent} \Psi(x \ot a) &=& x \ot {\mathfrak z}  a \\
\Psi(w \ot b) &=& w \ot {\mathfrak z}  b \nonumber \\
\Psi \big(\la \alpha, \be \ra \big) &=& \la {\mathfrak z} \alpha, 
\be\ra = \la \alpha, {\mathfrak z} \be \ra \nonumber
\end{eqnarray}
\noindent   for all $x \ot a \in \g \ot A$, $w \ot b \in W \ot B$, 
$\alpha, \be \in \mathfrak a$.
 \end{cor}
\psp

\begin{thm} \label{cent-rg} Let $\LL = ({\mathfrak g} \ot A)
\oplus (W \ot B) \oplus \langle \mathfrak a, \mathfrak a\rangle$ 
denote a Lie algebra graded by a finite reduced root system $\Delta$  
with coordinate algebra $\mathfrak a = A \oplus B$. Then $\cent(\LL) 
\cong \mathfrak Z_\mathfrak a$,  where $\mathfrak Z_\mathfrak a$ is 
the set of elements $\mathfrak z$ in $Z(\mathfrak a)\cap A$ which 
satisfy the following properties:
\begin{itemize}
 \item[{\rm(a)}]  $\la
{\mathfrak z}  \alpha, \be\ra = \la \alpha, {\mathfrak z} \be \ra$ 
for all $\al,\be \in \mathfrak a$;   
 \item[{\rm (b)}]
$\sum_t \la \al_t,\be_t \ra = 0$  implies  $\sum_t \la \mathfrak z 
\al_t,\be_t\ra  = 0$.
 \end{itemize}   More specifically, if $\mathfrak z
\in Z(\mathfrak a) \cap A$ satisfies (a) and (b),  and if  $\Psi \in 
\End_\F(\LL)$ is given by (\ref{cent})  above,  then $\Psi \in 
\cent(\LL)$;  and every element of $\cent(\LL)$ has this form.
\end{thm}

\proof   Suppose $\mathfrak z \in \mathfrak Z_\mathfrak a$, and 
define $\Psi$ as in (\ref {cent}). We need to know that the action 
of $\Psi$ on $\la \mathfrak a, \mathfrak a\ra$ is well-defined. But 
that is apparent from condition (b).   Also, to make sense of the 
definition, we must have  $\mathfrak z a \in A $ and $\mathfrak z b 
\in B$ for $a \in A$, $b \in B$, which is of course obvious in case 
$\mathfrak a = A$, i.e., $\De$ is simply laced. If $\De$ is not 
simply laced and $\De \ne  {\rm F}_4$ or ${\rm G}_2$,  the condition 
follows from the fact that the subspaces $A$ and $B$ are the 
symmetric and skew-symmetric elements with respect to an involution 
of $\mathfrak a$ (\cite{BZ}).  

In case $\De= {\rm F}_4$ or ${\rm G}_2$, the algebra $\mathfrak a$ 
is a unital algebra over the commutative associative subalgebra 
$A=A.1$, while the subspace $B$ is the kernel of an $A$-linear trace 
functional $\mathfrak a \to A$, and hence $AB\subset B$. The fact 
that $\Psi \in \cent(\LL)$ can be verified directly using 
Proposition~\ref{prop1rg}.  Now for the other direction, apply 
Corollary \ref{oneway} to deduce $\Psi$ has the form in 
(\ref{cent}).
  If  $\sum_t \la \al_t,\be_t \ra = 0$,
then $\Psi$ must map that expression to 0, so that  $\sum_t \la 
\mathfrak z \al_t,\be_t\ra  = 0$.   \qed  \psp

A  $\Delta$-graded Lie algebra $\LL$ with trivial centre has the 
following form  $\LL = ({\mathfrak g} \ot A) \oplus (W \ot B) \oplus 
D_{\mathfrak a,\mathfrak a}$, where $D_{\mathfrak a,\mathfrak a}$ is 
the space of all inner derivations. In this particular case,  
$D_{\mathfrak z \al, \be} = D_{\al, \mathfrak z \be}$ for all $\al, 
\be \in \mathfrak a$ and all $\mathfrak z \in Z(\mathfrak a) \cap 
A$.  Moreover, $\mathfrak zD_{\al, \be} = D_{ \mathfrak z \al,\be}$,  
so that (b) holds as well. Therefore, Theorem \ref{cent-rg} implies: 
\psp

\begin{cor} \label{cent-less}  Let $\mathfrak L = ({\mathfrak g} \ot A)
\oplus (W \ot B) \oplus D_{\mathfrak a,\mathfrak a}$ denote the 
centerless Lie algebra graded by a finite reduced root system 
$\Delta$ with coordinate algebra $\mathfrak a = A \oplus B$. Then 
$\cent(\LL) \cong Z(\mathfrak a) \cap A$.  \end{cor} \psp

Any  $\Delta$-graded Lie algebra $\KK = (\g \ot A) \oplus (W \ot B) 
\oplus \la \mathfrak a, \mathfrak a\ra$  with coordinate algebra 
$\mathfrak a = A \oplus B$ is a cover of the centreless 
$\Delta$-graded Lie algebra $\LL =   (\g \ot A) \oplus (W \ot B) 
\oplus D_{\mathfrak a, \mathfrak a}$ with that same coordinate 
algebra  via the map which is the identity on $(\g \ot A) \oplus (W 
\ot B)$ and sends $\la \al,\be \ra$ to $D_{\al,\be}$.    By Theorem 
\ref{cent-rg} and Corollary \ref{cent-less},   there is always an 
embedding

$$\cent(\KK)  \to Z(\mathfrak a) \cap A \cong
\cent(\LL).$$ \noindent Of course, we already knew that from Lemma 
\ref{oblem}(c). \psp
 
\begin{exam}\label{exaff} {\rm 
Let $\KK = \bigl( \g \ot \F[t,t^{-1}] \bigr) \oplus \F c$ be the 
derived algebra of an untwisted affine algebra. As an application of 
general results,  we have seen in Corollary \ref{remKM} and then 
again in Example \ref{examEALA}  that $\KK$ is central.   An 
alternate proof of this fact comes from specializing 
Theorem~\ref{cent-rg}.

Indeed, $\KK$  is a $\Delta$-graded Lie algebra with grading 
subalgebra the split simple Lie algebra $\g$.   The element $c$ is 
central, and $[x \ot t^m, y \ot t^n] = [x,y] \ot t^{m+n} + 
(x\,|\,y)\la t^m, t^n \ra$  for all $x,y \in \g$, where $\la t^m, 
t^n \ra =  m \delta_{m,-n}c$.  The centreless $\Delta$-graded Lie 
algebra $\LL$  with the same coordinate algebra  is just the loop 
algebra $\LL : =  \g \ot \F[t,t^{-1}] $, whose centroid according to  
Corollary~\ref{cent-less} is  $\cent(\LL) = \F[t,t^{-1}]$,  since 
$\mathfrak a = A =\F[t,t^{-1}]$ is a commutative associative algebra   
(compare Remark \ref{remloo} and Remark \ref{laurent} with $\sigma = 
\Id$).   By Theorem \ref{cent-rg}, we know that each element $\chi 
\in \cent(\KK)$ is determined by an element $\mathfrak z = \sum_p  
k_p t^p \in \mathfrak Z_{\mathfrak a} = \mathfrak a$, which 
satisfies    $\la\mathfrak z t^m, t^n \ra = \la t^m, \mathfrak z t^n 
\ra$ for all $m,n \in \Z$. Thus,

$$\sum_p k_p(m+p)\delta_{m+p,-n} = \sum_p k_p m \delta_{m,-n-p}.$$

\noindent must hold.  If $k_q \neq 0$ for some $q \neq 0$, then 
choosing $m,n$ so that $m=-n-q \neq 0$, we obtain $m+q = m$, a 
contradiction.  So it must be that $\mathfrak  z = k_0 1$. 
Consequently, $\cent(\KK) = \F \Id$.}   \end{exam}\pse

\subsection{Lie algebras graded by nonreduced root systems 
 (the BC$_r$ case)} \label{nonred} 
For simplicity,  we will assume $r \geq 3$.  Exceptional behavior 
occurs for small ranks, and somewhat different arguments need to be 
used for them.   We will not address those cases here.  In what 
follows,  we will apply results from \cite{ABG2} without 
specifically quoting chapter and verse.

\psp When $r \geq 3$, each BC$_r$-graded Lie algebra $\LL$ admits a 
decomposition,

$$\LL = (\g \ot A) \oplus (\mathfrak s  \ot B) \oplus (V \ot C) \oplus D,$$

\noindent relative to the grading subalgebra $\g$.   The spaces 
$\mathfrak s$ and $V$ are irreducible $\g$-modules and $D$ is the 
sum of the trivial $\g$-modules. Moreover,
\begin{itemize}
\item[{\rm (a)}]  The sum  $\mathfrak a = A \oplus B$ is a unital
algebra with involution $\eta$ whose symmetric elements are $A$ and 
skew-symmetric elements are $B$. \item[{\rm (b)}]   The algebra 
$\mathfrak a$ is associative in all cases except when $\g$ is of 
type C$_3$.  In that exceptional case $\mathfrak a$ is an 
alternative algebra, and the set $A$ of symmetric elements must be 
contained in the nucleus (associative center) of $\mathfrak a$. 
\item[{\rm (c)}] The space $C$ is a left $\mathfrak a$-module, and 
it is equipped with a hermitian or skew-hermitian form $\chi(\,,\,)$ 
depending on whether $\Delta_\g$ is of type B$_r$, D$_r$ or C$_r$. 
\item[{\rm (d)}]  $\mathfrak b = \mathfrak a \oplus C$ is an algebra 
(the {\em coordinate algebra} of $\LL$)  with product given by
$$(\al + c)\cdot (\al' + c') = \al\al' + \chi(c,c') + \al.c' +
{\al'}^\eta.c.$$
\end{itemize}

  Now suppose $\Psi \in \cent(\LL)$.   Then by Lemma
\ref{lemcr} (compare also Proposition \ref{prop1rg}),  it follows 
that there are transformations $\psi_{_\mathfrak a}$, $\psi_{_C}$, 
and $\psi_{_D}$ such that $$\Psi\bigl((x \ot a) + (s \ot b) + (v \ot 
c) + d\bigl) = (x \ot \psi_{_\mathfrak a}(a)) + (s \ot 
\psi_{_\mathfrak a}(b)) + (v \ot \psi_{_C}(c)) + \psi_{_D}(d)$$

Since $\mathfrak M =  (\g \ot A) \oplus (\mathfrak s  \ot B) \oplus 
\langle \mathfrak a,\mathfrak a\rangle$  is a subalgebra having 
exactly the same multiplication  as in (\ref{secondtype}) (think of 
$\mathfrak s$ as playing the role of $W$ in the reduced case), we 
obtain just as before that there exists an element $\mathfrak z \in 
\mathfrak Z_{\mathfrak a}$ such that $\Psi(x \ot a) = x \ot 
\mathfrak z a$, $\Psi(s \ot b) = s \ot \mathfrak z b$, and 
$\Psi(\langle \al, \al'\rangle) = \langle \mathfrak z \al, 
\al'\rangle = \langle \al, \mathfrak z \al' \rangle$ for all $a \in 
A$, $b \in B$, $\al,\al' \in \mathfrak a$, $x \in \g$, and $s \in 
\mathfrak s$.

Note that $ \Psi\bigl([x \ot 1, u \ot c] \bigr)   = \Psi\bigl(x.u 
\ot c\bigr) = x.u \ot \psi_{_C}(c)$, which must equal $[\Psi(x \ot 
1), u \ot c] = [x \ot \mathfrak z , u \ot c] = x.u \ot \mathfrak 
z.c$ for all $x \in \g$, $v \in V$, and $c \in C$. Thus, $\psi_C(c) 
= \mathfrak z.c$.  Applying the formulas in \cite[(2.8)]{ABG2}, we 
determine from the relation $\Psi\bigl([u \ot c, v \ot c']\bigr) = 
[\Psi(u \ot c), v\ot c'] = [u \ot c, \Psi(v \ot c')]$ that 
$\Psi\bigl(\langle c,c'\rangle \bigr) = \langle \mathfrak z.c, c' 
\rangle = \langle c, \mathfrak z.c'\rangle$ for all $c,c' \in C$.   
Since  $D = \langle \mathfrak a, \mathfrak a \rangle + \langle C, C 
\rangle$, this determines $\Psi$ completely.

 \psp
\begin{thm} \label{cent-nrg} Let $\LL = ({\mathfrak g} \ot A)
\oplus (\mathfrak s \ot B) \oplus (V \ot C) \oplus  \langle 
\mathfrak b, \mathfrak b\rangle$ denote a Lie algebra graded by 
$\hbox{\rm BC}_r$ for $r \geq 3$ with coordinate algebra $\mathfrak 
b = \mathfrak a \oplus C$ where $\mathfrak a = A \oplus B$.  Then 
$\cent(\LL) \cong \mathfrak Z_\mathfrak a$, where $\mathfrak 
Z_\mathfrak a$ is the set of elements in $Z(\mathfrak a)\cap A$ 
which satisfy the following properties:
\begin{itemize}
\item[{\rm(a)}]  $\la {\mathfrak z}  \be, \be'\ra = \la \be,
{\mathfrak z} \be' \ra$ for all $\be,\be' \in \mathfrak b$;
\item[{\rm (b)}]   $\sum_t \la \be_t,\be_t' \ra = 0$  implies
$\sum_t \la \mathfrak z \be_t,\be_t'\ra  = 0$ for all $\be_t,\be_t' 
\in \mathfrak b$.
 \end{itemize}   More specifically, if $\mathfrak z
\in Z(\mathfrak a) \cap A$ satisfies (a) and (b),  and if  $\Psi \in 
\End_\F(\LL)$ is given by
\begin{eqnarray}\label{centnr} \Psi(x \ot a) &=& x \ot {\mathfrak z}  a \\
\Psi(s \ot b) &=& s\ot {\mathfrak z}  b \nonumber \\
\Psi(v \ot c) &=& v \ot  {\mathfrak z}.c \nonumber \\
\Psi \big(\la \be, \be' \ra \big) &=& \la {\mathfrak z}  \be, 
\be'\ra = \la \be, {\mathfrak z} \be' \ra \nonumber
\end{eqnarray}
\noindent   for all $x \ot a \in \g \ot A$, $s \ot b \in \mathfrak s 
\ot B$, $v \ot c \in V \ot C$, and $\be, \be'  \in \mathfrak b$, 
then $\Psi \in \cent(\LL)$;  and every element of $\cent(\LL)$ has 
this form.
\end{thm}

\psp


\begin{thebibliography}{99}

 \bibitem[AABGP]{AABGP} B.N.~Allison, S.~Azam, S.~Berman, Y.~Gao,
 and A.~Pianzola, {\it Extended Affine Lie Algebras and Their
 Root Systems}, Mem. Amer. Math. Soc. \textbf {126}  vol.~603,
 Providence, R.I.  1997.

 \bibitem[ABG1]{ABG1}  B.N.~Allison, G.~Benkart, and Y.~Gao, Central
extensions of Lie algebras graded by finite root systems, Math. Ann. 
\textbf{ 316} ( 2000), 499-527.


\bibitem[ABG2]{ABG2}  B.N.~Allison, G.~Benkart, and Y.~Gao,
{\it Lie Algebras Graded by the Root Systems BC$_r$, $r \geq 2$}, 
Mem. Amer.  Math.  Soc.  \textbf {158} vol.~751, Providence, R.I. 
2002.

\bibitem[ABGP]{ABGP} B.N.~Allison, S.~Berman, Y. Gao and A. Pianzola,
A characterization of affine Kac-Moody Lie algebras, Comm. Math. 
Phys. \textbf{185} (1997), 671--688. 

\bibitem[ABP1]{ABP2} B.N.~Allison, S.~Berman and A.~Pianzola, Covering
algebras II: Isomorphism of loop algebras, J. Reine Angew. Math. 
\textbf{571} (2004), 39--71.

\bibitem[ABP2]{ABP3} B.N.~Allison, S.~Berman and A.~Pianzola, Iterated 
loop algebras, preprint February 2005.

\bibitem[BM]{BM} G.~Benkart and R.V.~Moody, Derivations, central
extensions, and affine Lie algebras, Algebras, Groups and Geom. 
\textbf{3} (1986)  456--492.

\bibitem[BO]{BO} G.~Benkart and J.M.~Osborn, Derivations and automorphisms
of nonassociative matrix algebras, Trans. Amer. Math. Soc. 
\textbf{263} (1981), 411--430.

\bibitem[BS]{BS} G.~Benkart and O.~Smirnov, Lie algebras 
graded by the root system BC$_1$,   J.  Lie Theory \textbf{13} 
(2003), 91--132.

\bibitem[BZ]{BZ}
G.~Benkart and E.~Zelmanov, Lie algebras graded by finite root 
systems and intersection matrix algebras, Invent. Math. \textbf{126} 
(1996), 1--45.



 \bibitem[BGKN]{BGKN}  S.~Berman, Y.~Gao,  Y.~Krylyuk, and
E. Neher, The alternative torus and the structure of elliptic 
quasi-simple Lie algebras of type $A_2$, Trans. Amer. Math. Soc. 
\textbf{347} (1995),  4315--4363.

\bibitem[BeM]{BeM} S.~Berman and R.V.~Moody,  Lie algebras
graded by finite root systems and the intersection matrix algebras 
of Slodowy, Invent. Math.  \textbf {108} (1992), 323--347.

\bibitem[Bo1]{bou:A} N.~Bourbaki, {\it Alg\`ebre}, Ch.~II,
Diffusion C.C.L.C, Paris 1970. 

\bibitem[Bo2]{bou:lieI} N.~Bourbaki, {\it Groupes et Alg\`ebres de Lie},
Ch.~I, Hermann, Paris 1971. 

\bibitem[Bo3]{bou:lieVI} N.~Bourbaki, {\it Groupes et Alg\`ebres de Lie},
Ch.~VI, Masson, Paris 1981.
 
\bibitem[F]{F} R.~Farnsteiner, Derivations and central extensions of finitely generated graded Lie algebras.  J. Algebra  \textbf {118}  (1988),  33--45.


\bibitem[J]{jake} N.~Jacobson, {\it Lie Algebras}, Interscience Publishers,
1962.

\bibitem[K]{K} V.G.~Kac, {\it Infinite Dimensional Lie Algebras},
Third Ed., Cambridge U. Press, Cambridge 1990.

\bibitem[L]{L} T.Y.~Lam, {\it A First Course in Noncommutative
Rings}, Graduate Texts in Math. vol. \textbf{131}, Springer-Verlag 
New York, 1991.

\bibitem[Mc]{Mc} K.~McCrimmon, {\it A Taste of Jordan Algebras},
Universitext, Springer-Verlag, New York 2004. 

\bibitem[Me]{Me} D.~Melville, Centroids of nilpotent Lie algebras,
Comm. Algebra \textbf{20(12)} (1992), 3649--3682.

\bibitem[MP]{MP} R.V.~Moody and A.~Pianzola, {\it Lie Algebras With
Triangular Decompositions}, Canad. Math. Soc. Series of Monographs 
and Advanced Texts, John Wiley \& Sons, New York 1995.

\bibitem[N1]{N1}   E.~Neher, Lie algebras graded by $3$-graded
root systems,  Amer.  J.~Math. \textbf{118} (1996),  439--491.

\bibitem[N2]{N2}  E.~Neher, Lie tori,  C.~R.~Math.~Acad.~Sci.~Soc.~R.~Can.  \textbf{26}  (2004),  no.~3, 84--89.
\bibitem[N3]{N3}  E.~Neher, Extended affine Lie algebras, C.~R.~Math.~Acad.~Sci.~Soc.~R.~Can. \textbf{26} (2004), no.~3,  90--96.

\bibitem[P]{P} D.~Passman, {\it Infinite Crossed Products},  
Pure and Applied Mathematics, \textbf{135}  Academic Press, Inc., 
Boston, MA, 1989.

\bibitem[Pi]{Pi} A.~Pianzola, Automorphisms of toroidal {L}ie algebras
and their central quotients, J. Algebra Appl. \textbf{1} (2002), 
113--121.


\bibitem[Po]{Po} K.~Ponomar{\"e}v, Invariant Lie algebras and Lie
algebras with a small centroid, Algebra Logic \textbf{40} (2001), 
365--377.

\bibitem[S1]{S1} K. Saito, Extended affine root systems 1 (Coxeter transformations),
Publ. RIMS, Kyoto Univ. \textbf{21} (1985) 75-179.

\bibitem[S2]{S2} K. Saito, Extended affine root systems 2 (flat invariants),
Publ. RIMS, Kyoto Univ. \textbf{26} (1990) 15-78.

\bibitem[Se]{Se} G.B. Seligman, {\it Rational Methods in  Lie Algebras},
Lect. Notes in Pure and Applied Math. \textbf {17} Marcel Dekker, 
New York 1976.

\bibitem[Sl]{Sl} P. Slodowy, A character approach to Looijenga's invariant
theory for generalized root systems, Compositio Math. \textbf{55} 
(1985), 3-32. 

\bibitem[T]{T} J.~Tits, Une classe d'alg\'ebres
de Lie en relation avec les alg\'ebres de Jordan, Indag. Math. 
\textbf{24} (1962), 530--535.

\bibitem[Y1]{Y2} Y.~Yoshii, Root-graded Lie algebras with compatible
grading, Comm. Algebra \textbf{29} (2001), 3365--3391.

\bibitem[Y2]{Y3} Y.~Yoshii, Classification of division
{$\mathbb Z\sp n$}-graded alternative algebras, J. Algebra 
\textbf{256} (2002), 28--50.

\bibitem[Y3]{Y4} Y.~Yoshii, Root systems extended by an abelian group
and their Lie algebras,  J. Lie Theory  \textbf{14} (2004), 
371--394.

\bibitem[Y4]{Y5} Y.~Yoshii, Lie tori -- A simple characterization of 
extended affine Lie algebras, preprint 2003. 
\end{thebibliography}
\end{document}